\documentclass{gtpart}     
\usepackage{pinlabel}  

\usepackage{amsfonts}
\usepackage{mathrsfs} 
\usepackage{enumerate}
\usepackage{type1cm}

\title{Geodesible contact structures on 3--manifolds}

\author{Patrick Massot}
\givenname{Patrick}
\surname{Massot}
\address{\'Ecole Normale Sup\'erieure de Lyon, 69364 LYON Cedex 07, France}
\email{patrick.massot@umpa.ens-lyon.fr}
\urladdr{http://www.umpa.ens-lyon.fr/~pmassot/}

\keyword{contact structures}
\keyword{totally geodesic}
\keyword{Seifert manifolds}
\keyword{twisting number}
\subject{primary}{msc2000}{57M50}
\subject{secondary}{msc2000}{57R17}

\arxivreference{0711.0377}
\arxivpassword{fuf9n}

\volumenumber{}
\issuenumber{}
\publicationyear{}
\papernumber{}
\startpage{}
\endpage{}
\doi{}
\MR{}
\Zbl{}
\received{}
\revised{}
\accepted{}
\published{}
\publishedonline{}
\proposed{}
\seconded{}
\corresponding{}
\editor{}
\version{}

\newtheorem*{cor-intro}{Corollary}    
\newtheorem{thm}{Theorem}[section]    
\newtheorem{mythm}{Theorem}
\newtheorem{lem}[thm]{Lemma}          
\newtheorem{prop}[thm]{Proposition}   
\newtheorem{cor}[thm]{Corollary}      
\theoremstyle{definition}
\newtheorem{defn}[thm]{Definition}    
                                    
\newtheorem{rem}[thm]{Remark}         
\newtheorem{ex}[thm]{Example}

\def\N{\mathbb{N}}  
\def\Z{\mathbb{Z}}  
\def\R{\mathbb{R}}  
\def\C{\mathbb{C}}  
\def\Q{\mathbb{Q}}  
\def\Card{\mathrm{Card}} 
\def\tight{\mathrm{Tight}} 
\def\F{\mathscr{F}}  
\def\tor{\mathrm{Tor}}
\def\SL{\mathrm{SL}}
\def\tr{\mathrm{tr}}
\def\dprime{d\hspace{.5pt}'}
\def\dhat{\skew7\hat{d}}
\def\gammabar{\skew2\bar{\gamma}}
\def\xihat{\skew2\hat{\xi}}
\def\xitilde{\skew4\tilde{\xi}}
\def\Vtilde{\skew3\widetilde{V}}
\def\Vhat{\skew{1.5}\widehat{V}}
\def\Lhat{\skew{1.5}\widehat{L}}
\def\Wtilde{\skew3\widetilde{W}}
\def\phibar{\skew4\bar{\phi}}
\def\STB{ST^*\kern-1pt B}

\def\P{\mathcal{P}}  
\def\S{\mathrm{S}}  

\newcommand{\step}[1]{{\flushleft\textsc{Step #1}}}

\begin{document}

\begin{abstract}    
	In this paper, we study and almost completely classify contact structures
on closed 3--manifolds which are totally geodesic for some Riemannian metric.
Due to previously known results, this amounts to classifying contact structures 
on Seifert manifolds which are transverse to the fibers. Actually, we obtain the
complete classification of contact structures with negative (maximal) twisting
number (which includes the transverse ones) on Seifert manifolds whose base is 
not a sphere, as well as partial results in the spherical case.
\end{abstract}

\maketitle


\section*{Introduction}

A geodesible plane field on a manifold is a plane field $\xi$ which is totally 
geodesic for some Riemannian metric $g$, which means that every
geodesic of $g$ that is somewhere tangent to $\xi$ is everywhere
tangent to $\xi$. If such a plane field is integrable then
(locally) its leaves are totally geodesic submanifolds for $g$
and this class of foliations has been much studied. Following a
suggestion of \'E~Ghys, we investigate here, by contrast,
geodesible contact structures on 3--manifolds. The basic example
of a totally geodesic contact structure is the standard contact
structure on the sphere which, by definition, is the plane field
orthogonal to the Hopf fibers for the canonical Riemannian
metric. The goal of this paper is to give an almost complete
classification of geodesible contact structures on closed
3--manifolds.

The existence of a geodesible plane field is a strong constraint on the topology
of the underlying manifold. In dimension $3$ for instance, a theorem due to Yves
Carri\`ere shows that, up to diffeomorphism, a closed manifold equipped with a 
coorientable geodesible plane field is either a Seifert manifold with a plane 
field transverse to its fibers or a torus bundle over the circle whose linear 
monodromy $A \in \SL_2(\Z)$ satisfies $\tr(A) > 2$, this manifold being endowed
with a plane field transverse to the foliation spanned by an
eigendirection of $A$ (which is the strong stable or unstable foliation of the
Anosov flow given by $A$), see Section \ref{S:plane_fields}. On torus bundles, 
contact structures are completely classified by E~Giroux and K~Honda 
in \cite{Giroux_2000, Honda_II}. This 
classification implies that a geodesible contact structure exists if and only if
$\tr(A) > 2$ (as for any coorientable plane field) and that, in this case, it is
unique up to isotopy. On Seifert manifolds, in contrast, the situation is much 
more delicate: currently there is no general classification of
tight contact structures and our main task will be to classify
those which are transverse to the fibers---since these are the
geodesible ones, up to isotopy. In order to state the
results, we define now the twisting number of a contact structure
on a Seifert manifold.

Let $V$ be a Seifert 3--manifold and $K \subset V$ a regular fiber, that is, a 
fiber admitting a trivialized neighborhood $D^2 \times \S^1 \supset \{0\} \times 
\S^1 = K$ in which all circles $\{\cdot\} \times \S^1$ are also Seifert fibers. Then
$K$ has a canonical (homotopy class of) normal framing given by this splitting.
Now let $\xi$ be a contact structure on $V$ and $L$ a Legendrian curve smoothly
isotopic to $K$. Given a smooth isotopy $\varphi$ taking $K$ to $L$, let $t (L,
\varphi)$ be the difference between the contact framing of $L$ and the image by
$\varphi$ of the canonical framing of $K$. Then define the twisting number of
$L$ by $t(L) = \sup_\varphi t(L,\varphi)$ and the twisting number of $\xi$ to be
the supremum of the twisting numbers $t(L)$ for all Legendrian curves isotopic 
to regular fibers. One can prove that overtwisted contact structures have
infinite twisting number.

\begin{mythm}[Sections \ref{S:debut_charac} and \ref{S:charac}]
  \label{thm:charac}
Let $V$ be a closed Seifert 3--manifold. A contact structure $\xi$ on $V$ is 
isotopic to a transverse one 
if and only if it is universally tight and has negative twisting number. On the
other hand, any transverse contact structure $\xi$ on $V$ is symplectically
fillable. 
\end{mythm}

It turns out that, when it exists, the (unique) geodesible contact structure on 
a torus bundle is also universally tight and symplectically fillable. Hence
geodesibility is (apparently the first example of) a compatibility condition 
between a Riemannian metric and a contact structure which implies tightness (and
symplectic fillability).

Note here that all our Seifert manifolds are oriented with oriented fibers, and
that all contact structures we consider are positive for this orientation and
(co)\,oriented unless explicitly stated otherwise.

The next step is to determine which Seifert manifolds admit transverse contact 
structures. For this we need a concrete description of Seifert manifolds. Start
with a compact oriented surface $R$ of genus $g$ with $r+1$ boundary components, 
take $R \times \S^1$ with the product orientation and attach $r+1$ solid tori 
$W_0, \dots, W_r$ to its boundary, in order to get a closed manifold. The gluing 
is prescribed by integers $b, (\alpha_1, \beta_1), \dots, (\alpha_r, \beta_r)$ 
as follows: the meridian disk of $W_i$ is attached to a curve
whose homology class is $-\alpha_i [\partial_iR] + \beta_i
[\S^1]$, $0 \le i \le r$, where $(\alpha_0, \beta_0) = (1,b)$,
$1\le \beta_i<\alpha_i$ for $i>0$ and $\partial_i R$ is the
$i$-th boundary component of $R$ with the induced orientation
(see Subsection~\ref{subsection:notations}).
The Seifert invariants of the resulting manifold are, by definition, $(g, b,
(\alpha_1, \beta_1), \dots, (\alpha_r, \beta_r))$.

\begin{mythm}[Section \ref{S:existence}]
  \label{thm:existence}
Let $V$ be a Seifert manifold with invariants
\[ (g, b, (\alpha_1, \beta_1), \dots, (\alpha_r, \beta_r)) \]
and let $n$ be a positive integer. The following properties are equivalent:
\begin{itemize}
\item[(i)] 
$V$ carries a transverse contact structure whose twisting number is at least 
$-n$;
 \item[(ii)] 
$V$ carries a contact structure whose twisting number is negative and is at 
least $-n$;
\item[(iii)] 
there exist integers $x_0, \dots, x_r$ such that $\sum x_i = 2-2g$ and
$\frac{x_i-1}{n} < \frac{\beta_i}{\alpha_i}$ for $0 \le i \le r$, with $\alpha_0 = 1$ and 
$\beta_0 = b$.
\end{itemize}
\end{mythm}

Using the equivalence of (i) and (iii), one can recover the
following criterion:

\begin{cor-intro}
A Seifert manifold with invariants
\[ (g, b, (\alpha_1, \beta_1), \dots, (\alpha_r, \beta_r)) \]
carries a transverse contact structure if and only if one of the following 
conditions holds:
\begin{itemize}
\item 
$-b-r \leq 2g-2$;
\item 
$g = 0$, $r \leq 2$ and $-b-\sum \beta_i/\alpha_i < 0$;
\item 
$g=0$, $r \geq 3$, $-b-r=-1$ and there exist relatively prime
integers $0< a < m$ and a reordering of the pairs $(\alpha_i,\beta_i)$, 
$1 \le i \le r$, such that:
\[ \frac{\beta_1}{\alpha_1} > \frac{m-a}{m}, \quad 
   \frac{\beta_2}{\alpha_2} > \frac{a}{m}, \quad \text{and} \quad
   \frac{\beta_i}{\alpha_i} > \frac{m-1}{m} \ \forall i \ge 3. \]
\end{itemize}
\end{cor-intro}

This corollary was first proved for circle bundles ($r=0$) by 
E~Giroux in \cite{Giroux_2001} and A~Sato and T~Tsuboi in
\cite{Tsuboi_94}. Then the case of general Seifert manifolds was
treated by K~Honda in \cite{Honda_98} (almost) and by P~Lisca
and G~Mati\'c in \cite{Lisca_Matic}. The proof in
\cite{Honda_98} is direct and relies on subtle results on circle
diffeomorphisms established in \cite{EHN, JN_1, JN_2, Naimi}. The
proof in \cite{Lisca_Matic} combines the Eliashberg--Thurston
perturbation theorem with the more analytical 
adjunction inequality in symplectic geometry. In the case $g=0$,
$r=3$ and $b=-2$, a topological proof of
Theorem~\ref{thm:existence} by P~Ghiggini appeared in
\cite{Paolo_negative} while this paper was in preparation. 

The main tool in our work is a set of normal forms which we now
describe. Take a Seifert manifold $V = (R \times \S^1) \cup W_0
\cup \dots \cup W_r$ and fix a complex structure $J$ on $R$ which
defines the orientation of $R$. For each non-singular 1--form
$\lambda$ on $R$, the Pfaff equation 
\[ \cos(n\theta) \, \lambda + \sin(n\theta) \, \lambda \circ J = 0 \]
defines a contact structure on $R \times \S^1$ and any contact structure $\xi$ on
$V$ which extends this one will be denoted by
\[ \xi = \xi (\lambda, n, \xi_0, \dots, \xi_r) \qquad
   \text{where} \quad \xi_i = \xi_{|W_i}, \quad 0 \le i \le r \,. \]
If, in addition, the integer $-n$ is equal to the twisting number $t(\xi)$ then
we say that $\xi$ is in normal form.

The starting point of our study is the observation that any contact structure on
$V$ with negative twisting number is isotopic to one which is in normal form. In
particular, the integers $x_i$ in condition (iii) of the existence criterion are
the indices of $\lambda$ along the boundary components of $R$. The collection of
these indices is called the multi-index of the normal form.

\medskip

Normal forms and their multi-indices will be our main tools to classify contact
structures with negative twisting numbers. The following theorem determines the 
possible negative twisting numbers for a contact structure on a given Seifert 
manifold provided the base surface has positive genus. 

\begin{mythm}[Section \ref{S:calcul_e}]
\label{thm:e}
Let $V$ be a Seifert manifold with invariants
\[ (g, b, (\alpha_1, \beta_1), \dots, (\alpha_r, \beta_r)), \quad g>0, \]
and let $\xi = \xi (\lambda, n, \xi_0, \dots, \xi_r)$ be a contact structure
on $V$ with multi-index $(x_0, \dots, x_r)$. Then $\xi$ is 
in normal form---that is, $t(\xi)=-n$---if and only if, for $0 \le i \le n$,
the following conditions hold:
\begin{itemize}
\item
the contact structure $\xi_i = \xi_{|W_i}$ is tight;
\item $(x_i-1)/n<\beta_i/\alpha_i$ and
the triangle in $\R^2$ with vertices $(0,0)$, $(\alpha_i,\beta_i)$, $(n,x_i-1)$ 
does not contain any integer point whose abscissa is less than
$n$ except possibly its vertices---here again, $\alpha_0 = 1$ and $\beta_0 = b$.
\end{itemize}
\end{mythm}

Note that the second condition is automatically fulfilled when $n=1$ and can be
also expressed in arithmetic terms using the continued fraction
expansions of the $\beta_i/\alpha_i$'s, see Subsection~\ref{SS:local_e}. 
When $g=0$, $t(\xi)=-n$ implies the above
conditions but the converse is not known.

To classify contact structures with negative twisting number on a given Seifert 
manifold, it remains to understand when two contact structures in normal form 
with the same twisting number are isotopic. Here two cases appear
which require different approaches.

Recall that a tight contact structure $\xi$ on a solid torus $W$ whose boundary
is $\xi$--convex has a relative Euler class in $H^2 (W, \partial W)$: it is the 
obstruction to extending inside $W$ a non-singular vector field on $\partial W$ 
which is tangent to $\xi_{|\partial W}$---tightness implies the existence of 
such a vector field. 

\begin{mythm}[Section \ref{S:classification}]
\label{thm:flexible}
Let $V$ be a Seifert manifold with invariants
\[ (g, b, (\alpha_1,\beta_1), \dots, (\alpha_r,\beta_r)) \]
and assume that $-b-r < 2g-2$.

Every contact structure on $V$ with twisting number $-1$ has a
normal form whose multi-index is $(2-2g-r, 1, \dots, 1)$ and has
a Stein filling whose underlying smooth manifold depends only on
$V$. 

Two contact structures $\xi (\lambda, 1, \xi_0, \dots, \xi_r)$
and $\xi (\lambda', 1, \xi_0', \dots, \xi_r')$ in normal form 
with the above multi-index are
isotopic if and only if each $\xi_i$ has the same relative Euler
class as $\xi_i'$.

Among all isotopy classes of contact structures with twisting
number $-1$, exactly two contain transverse contact
structures---and only one if we consider non-oriented contact
structures.
\end{mythm}

In the $g=0$ case, tight contact structures on these manifolds were
classified by H~Wu in \cite{Wu_classif} (using a similar strategy
to distinguish isotopy classes). The case of circle bundles
($r=0$) was also previously treated by E~Giroux in
\cite{Giroux_2001} but without the precise counting of isotopy
classes and by K~Honda in \cite{Honda_II} but with some gap in
this counting. To prove here the second part of the theorem, we
use a result of P~Lisca and G~Mati\'c \cite{LM_Stein}: the
contact structures are distinguished by the first Chern classes
of their Stein fillings.

In the next theorem, we denote by $\lceil x \rceil$ the smallest integer which 
is not less than $x$. We also call $R$--class of a contact structure in normal
form $\xi (\lambda, n, \xi_0, \dots, \xi_r)$ the homotopy class of $\lambda$
among non-singular 1--forms on $R$.

\begin{mythm}[Section~\ref{S:classification}]
\label{thm:rigide}
Let $V$ be a Seifert manifold with invariants
\[ (g, b, (\alpha_1,\beta_1), \dots, (\alpha_r,\beta_r)), \qquad g>0, \]
let $n$ be a positive integer and assume that $-b-r = 2g-2$ or $n>1$. 

Every contact structure on $V$ with twisting number $-n$ has a normal form 
and any such normal form has multi-index 
$(nb, \lceil n\beta_1/\alpha_1 \rceil, \dots, \lceil n\beta_r/\alpha_r \rceil)$.

Two contact structures $\xi (\lambda,  n, \xi_0, \dots, \xi_r)$
and $\xi (\lambda',  n, \xi_0' ,\dots, \xi_r')$ in normal form
are isotopic if and only if they have the same $R$--class and each
$\xi_i$ has the same relative Euler class as $\xi_i'$.

Among all isotopy classes of contact structures with twisting
number $-n$ and a fixed $R$--class, exactly one if $n\beta_i\equiv
1 \mod \alpha_i$ for all $i$ and two otherwise contain transverse
contact structures---and only one in all cases if we consider
non-oriented contact structures.
\end{mythm}

The case of circle bundles ($r=0$) was treated by E~Giroux and
K~Honda in \cite{Giroux_2001, Honda_II}. On the other hand, in
the case $g=1$ and $r=1$, tight contact structures on $V$  were
classified by P~Ghiggini in \cite{Ghiggini_thesis} but without the twisting
number computation and without determining which contact
structures are universally tight or transverse.

A key step in the proof of the above theorem is to show that,
given a contact structure with twisting number $-n$ under our
hypotheses, two Legendrian curves which are smoothly isotopic to
the regular fibers and have twisting number $-n$ are Legendrian
isotopic. An analogous statement was obtained by J~Etnyre and 
K~Honda in \cite{Etnyre_Honda_knots} for Seifert structures on
$\S^3$---but the result was a corollary of the
classification---and by P~Ghiggini in \cite{Ghiggini_T3} for
Seifert structures on $T^3$.

According to our existence criterion, the two theorems above classify all 
contact structures with negative twisting number on any Seifert manifold whose 
base is a  surface of positive genus.

Throughout the paper, we assume the reader is familiar with the
theory of $\xi$--convex surfaces developed by E~Giroux in
\cite{Giroux_91} but we will briefly recall in Sections
\ref{S:torus_bundles} and \ref{S:toric_quotients} the results we
need about the classification of tight contact structures on
toric annuli and solid tori.

\paragraph{Acknowledgements} We warmly thank Emmanuel Giroux for
suggesting the topic of this paper, for many interesting
discussions and for his constant encouragement. We also thank Ko
Honda and Paolo Ghiggini for interesting discussions related to
this paper and the referee for pointing out several small
inaccuracies and helping improve the exposition. This work was
partially supported by the ANR Symplexe.

\section{Geodesible plane fields}
\label{S:plane_fields}

Totally geodesic contact structures were considered by R~Lutz and
T~Hangan for constant curvature Riemannian metrics in
\cite{Lutz_Hangan} and by R~Lutz alone with extra relations
between a contact structure and a arbitrary metric in
\cite{Lutz_helico}. Here by contrast we consider contact
structures which are totally geodesic for some arbitrary 
Riemannian metric with no extra condition. Besides the theory of
$\xi$--convex surfaces---which was not available to Lutz---what
makes this generality possible is the relation between geodesible
plane fields and a certain class of 1--dimensional foliations. 
This relation gives the following topological characterization 
due to Y~Carri\`ere \cite{Carriere}.

\begin{prop}
\label{prop:carac_topo}
  A coorientable plane field $\xi$ on a closed orientable
3--manifold $V$ is geodesible if and only if the pair
$(V,\xi)$ is isomorphic to one of the following:
\begin{itemize}
  \item a Seifert manifold endowed with a plane field transverse to its fibers;
  \item a hyperbolic torus bundle $T^3_A$ with monodromy $A$---where 
	$A\in SL_2(\Z)$ and $tr(A) > 2$---equipped with a plane field
	transverse to the foliation spanned by one of the
	eigendirections of $A$. 
\end{itemize}
\end{prop}

\noindent Given a
geodesible plane field $\xi$ this relation also explains how the
metrics $g$ such that $\xi$ is totally geodesic for $g$ look
like. In the last subsection we define our notations about
Seifert manifolds.

\subsection{Totally geodesic plane fields and Riemannian foliations}

Throughout this section, $\F$ will be a 1--dimensional foliation on a 3--manifold
and $F$ the corresponding line field.

\begin{defn}
  A codimension 2 foliation is \emph{Riemannian} 
 if it admits transverse disks equipped with Riemannian metrics such that
 each leaf meets at least one disk and holonomy maps are isometries.
\end{defn}

This definition and the following proposition essentially go
back to \cite{Reinhart}. They can be generalized to arbitrary
dimensions and codimensions.

\begin{prop}
\label{prop:geod_riem}
  A plane field on a 3--manifold is geodesible if and only if
it is transverse to some 1--dimensional Riemannian foliation.
\end{prop}

\begin{proof}
First suppose that a plane field $\xi$ on a 3--manifold $M$ is
transverse to a Riemannian foliation $\F$\!. 
We equip $\xi$ with the Riemannian metric $g$ pulled back from
the transverse disks by holonomy and then extend this metric to
$TM$ such that $\xi$ is orthogonal to $F$ and $F$ is given an
arbitrary metric.
We now show that $\xi$ is totally geodesic for $g$.
  
  It is sufficient to prove that for every $(x,v)\in \xi$ there
exists a small curve $\gamma$ tangent to $\xi$,
distance minimizing and satisfying $\gamma(0)=x$ and
$\gamma'(0)=v$. Let $T$ be a transverse disk intersecting the
leaf through $x$ and $\pi$ the projection from a neighborhood of
$x$ to $T$ along the leaves. Let $\gammabar$ be a small geodesic
in $T$ starting at $\pi(x)$ with velocity vector $\pi_*v$ and
$\gamma$ be the unique lift of $\gammabar$ to $M$ which is
tangent to $\xi$ and set $y=\gamma(t)$ for some small $t$. For
every curve $\alpha$ with $\alpha(0)=x$ and $\alpha(t)=y$ we
denote by $\alpha'_\xi$ the orthogonal projection of $\alpha'$ onto
$\xi$. 
The curve $\gamma$ is not longer  than $\alpha$ because
\[\int |\alpha'| \geq \int |\alpha'_\xi|=
\int |(\pi\circ\alpha)'|=l(\pi\circ \alpha)\geq l(\gammabar)=
\int|\gammabar'|=\int|\gamma'|
\]
where $l$ is the length function on $T$ and we used that 
$\pi_*:\xi\to TT$ is isometric by construction and $\gammabar$ is
(locally) minimizing.

Conversely, suppose that $\xi$ is totally geodesic for some metric
$g$. We show that $\xi^\perp$ integrates to a Riemannian foliation
$\F$.

We take any system of transverse disks for $\F$ intersecting all
leaves. Each point $x$ on a
leaf intersecting a transverse disk $T$ at a point $y$ defines a
metric on $T_yT$ by pushing $g|_\xi$ using the infinitesimal
holonomy. We have to show that this metric is independent of $x$
on a given leaf. It is sufficient to prove that if $x$ and $x'$
are on the same plaque $L$ then the infinitesimal holonomy is an
isometry from $\xi_x$ to $\xi_{x'}$.

The geodesics orthogonal to the plaques are tangent to $\xi$ so all the
nearby plaques are contained in tubes around $L$ according to the
(generalized) Gauss lemma. This implies
that $\F$ can be locally parametrized by a flow which preserves
the tubes and the fibration in disks of a tubular neighborhood of
$L$. This flow then preserves the distance to $L$ and so the
infinitesimal holonomy from $\xi_x$ to $\xi_{x'}$ is an isometry.
\end{proof}

\subsection{The Carri\`ere classification}

  If a 1--dimensional foliation is directed by a Killing vector
field for some Riemannian metric $g$ then the induced metrics 
on a suitable system of transverse disks intersecting all leaves are
invariant by holonomy so this foliation is Riemannian.
Such a foliation is called isometric.

\begin{ex}
\label{example:seifert}
Any foliation given by a locally free action of $\S^1$ is Riemannian
since any Riemannian metric can be averaged to give an invariant
metric.
\end{ex}

This gives a canonical Riemannian foliation on every Seifert
manifold where the leaves are the fibers. Amongst such manifolds,
lens spaces (including $\S^3$ and $\S^2\times \S^1$) also have
deformed versions of the preceding example which are still
isometric foliations:

\begin{ex}
\label{example:lens_spaces}
  We view $\S^3$ sitting in $\C^2$ as 
  \[\S^3=\bigl\{(z_1,z_2)\in \C^2\;;\; |z_1|^2+|z_2|^2=1\bigr\}.\]
For every pair $(\lambda,\mu)$ of real numbers there is a flow
$\phi_t(z_1,z_2)=(e^{i\lambda t}z_1,e^{i\mu t}z_2)$
which defines a Riemannian foliation. If $\lambda$ is a rational
multiple of $\mu$ then we get an example of the preceding type.
The lens space $L_{p,q}$ is the quotient of $\S^3$ by 
$(z_1,z_2)\mapsto (e^{2i\pi/p}z_1,e^{2i\pi q/p}z_2)$ and $\phi_t$
clearly induces a Riemannian foliation on the quotient. These
foliations are always $C^\infty$ close to a Seifert foliation.
\end{ex}

The deformed example on $\S^2\times \S^1$ is:

\begin{ex}
  \label{example:s2s1}
  Let $f$ be a rotation of $\S^2$ with irrational angle. The
foliation of $\S^2\times\R$ by $\{\cdot\}\times\R$ induces an
isometric foliation of 
\[
\frac{\S^2\times\R}{(x,t)\sim(f(x),t+1)}\simeq \S^2\times \S^1
\]
which is $\mathcal C^\infty$--close to a Seifert fibration.
\end{ex}

The following examples are more exotic since they are
not isometric foliations.

\begin{ex}
 \label{example:t3a}
 For every matrix $A\in \SL_2(\Z)$, if $\tr(A)>2$ then $A$ has
 two positive eigenvalues $\lambda,1/\lambda$. The eigendirection
 corresponding to $\lambda$ defines a foliation on $T^2\times\R$
 which induces a 1--dimensional Riemannian foliation on $T^3_A$,
 see \cite{Carriere}. Remark that if $\tr(A)<-2$ then one gets a
 non-orientable foliation so we are not interested in this case
 because we only consider coorientable contact structures.
\end{ex}

\begin{proof}[Proof of Proposition~\ref{prop:carac_topo}]
  According to Proposition~\ref{prop:geod_riem}, a plane field is
geodesible if and only if it is transverse to a Riemannian
foliation. According to Carri\`ere's classification in \cite{Carriere}, 
any closed oriented 3--manifold equipped with a 1--dimensional Riemannian
foliation is diffeomorphic to one of the following:
\begin{itemize}
\item the torus $T^3$ with a linear foliation,
\item a lens space with a foliation of Example~\ref{example:lens_spaces},
\item the product $\S^2 \times \S^1$ with a foliation of Example~\ref{example:s2s1},
\item a Seifert manifold with its fibration, see Example~\ref{example:seifert},
\item a torus bundle $T^3_A$, $\tr(A)>2$, with the foliation of
  Example~\ref{example:t3a}.
\end{itemize}
Because transversality is an open condition, a
plane field is transverse to one of the first three types of
foliations if and only if it is transverse to a foliation of the
fourth type.
\end{proof}

\subsection{Seifert manifolds}
\label{subsection:notations}
We recall the definitions of Seifert manifolds and their
invariants to fix conventions for notations and orientations. We
then consider two classes of examples: the Seifert structures of
$\S^3$ and the bundles of cooriented contact elements of
2--dimensional orbifolds.

Conventions for Seifert invariants vary greatly from papers to
papers. For instance, here we use the same one as P~Lisca and G~Mati\'c
in \cite{Lisca_Matic} but H~Wu uses a different one 
in \cite{Wu_classif}. To go from one convention
to another it suffices to use the
normalization $1\leq \beta_i <\alpha_i$
and the fact that the rational and
integer Euler numbers---denoted by $e$ and $e_0$ respectively in
the following---are the same with any convention. In particular,
$\beta_i/\alpha_i$ in one convention can become
$1-\beta_i/\alpha_i$ in another one but this can be checked using $e$ and
$e_0$. In the notation $M(r_1,\dots,r_k)$ used for instance in
\cite{Wu_classif}, $e=\sum r_i$.
A useful exercise in conversion is to compare the
conventions in the  discussion following the proof of 
Theorem~\ref{thm:e} at the end of Section~\ref{S:calcul_e} to the
conventions in the cited papers. 

Let $B$ be a closed oriented 2--orbifold of
genus $g$ with $r$ elliptic points $f_1,\dots,f_r$ of order
$\alpha_1,\dots,\alpha_r$ and no other exceptional point. Let
$f_0$ be a regular point of $B$ and $D_0,\dots,D_r$ be 
pairwise disjoint closed disks such that $f_i$ is in $D_i$ for every $i$.
Denote by $R$ the smooth surface $B\setminus \bigcup \mathring D_i$. In the
following, boundary components of $R$ will always be oriented as
boundary components of the $D_i$'s.

Let $V'$ be $R\times \S^1$, fix an orientation for the $\S^1$
factor and use the product orientation on $V'$.
The first homology group of the boundary components of $V'$ have a
basis $(S_i,F_i)$ where $S_i$ is the homology class of
$\partial D_i\times\{\cdot\}$ and $F_i$ is the homology class of
$\{\cdot\}\times \S^1$ ($S$ stands for section and $F$ for fiber). We
orient these boundary components by imposing the intersection
number $S_i\cdot F_i=+1$.

The Seifert manifold $V$ with invariants
$(g,b,(\alpha_1,\beta_1),\dots,(\alpha_r,\beta_r))$ is, by
definition, the manifold obtained by filling each boundary torus
of $V'$ with a solid torus $W_i$ which has a meridian whose
homology class is
$M_i=\alpha_i S_i+\beta_i F_i$ with
$(\alpha_0,\beta_0)=(1,b)$ and $0<\beta_i<\alpha_i$ for $i>0$.

Each solid torus $W_i$ is identified with $W(\alpha_i,\beta_i)$
where $W(\alpha,\beta)$ is the quotient of
$\R/\Z\times\R/\Z\times[0,1]$ by the relation 
$(x,y,t)\sim(x',y',t')$ if $t=t'=0$ and $\alpha(y'-y)=\beta(x'-x)$.
The solid torus $W(\alpha,\beta)$ is foliated by the images of
the circles $\{\cdot\}\times\R/\Z\times\{\cdot\}$. This foliation is
called the standard Seifert fibration on $W(\alpha,\beta)$.

These foliations inside the $W_i$'s  extend the foliation of $V'$ by $\S^1$ 
to a foliation of $V$ such that $B$ is the orbit space. The leaves
of this foliation are called the fibers of $V$. Smooth points of
$B$ are called regular fibers whereas elliptic points are called
exceptional fibers. A closed subset of a Seifert manifold is said
to be vertical if it is isotopic to a union of regular fibers.

The rational Euler number of $V$ is 
$e(V)=-\sum_{i=0}^r\frac{\beta_i}{\alpha_i}$
and the integer Euler number is $e_0(V)=-b-r$.

\begin{ex}[Seifert structures on $\S^3$]
\label{ex:s3}
 Let $(L_0,L_1)$ be two fibers of the Hopf fibration. One can
see $\S^3$ as a quotient of $T^2\times[0,1]$ such that the
projection map $p$ is a diffeomorphism from $T^2\times(0,1)$ to
$\S^3\setminus (L_0\cup L_1)$ and collapses $T_0$ to $L_0$ and
$T_1$ to $L_1$---where we set $T_t=T^2\times\{t\}$. For each
point $x\in L_0$ (resp. $x\in L_1$) $p^{-1}(x)$ is a circle
$\{pt\}\times \R/\Z\subset T_0$ (resp.
$\R/\Z\times\{pt\}\subset T_1$). Any foliation of $T^2\times
(0,1)$ by circles not homologous to the above circles project
to a foliation of $\S^3\setminus(L_0\cup L_1)$ which, together
with $L_0$ and $L_1$ gives a Seifert structure on $\S^3$. These
are the only Seifert structure on $\S^3$, see eg \cite{Scott}.
\end{ex}

\begin{ex}[Cooriented contact elements bundles]
Seifert manifolds are the total spaces of circle bundles over
2--orbifolds in the sense of Thurston. One type of such bundles
is especially interesting for us, the cooriented contact elements
bundle $\STB$ of an orbifold $B$ where we orient the fibers
unlike the natural orientation.  It has Seifert invariants 
\[(g,2-2g-r ,(\alpha_1,1),\dots,(\alpha_r,1))\]
and rational Euler number 
$e(\STB)=-\chi(B)$ where the orbifold Euler number of $B$ is
$$\chi(B)=\chi(R)+\sum_{i=0}^r\frac{1}{\alpha_i}
=2-2g-r+\sum_{i=1}^r\frac{1}{\alpha_i}.$$

As in the special case of smooth
surfaces, $\STB$ carries a canonical (positive) contact structure
denoted by $\xi_B$.
\end{ex}

\section{Contact structures on toric annuli}
\label{S:torus_bundles}

Here we recall some facts about tight contact
structures on toric annuli which are needed in this text. 
All the results here are contained in \cite{Giroux_2000} or
follow directly from results therein (see also \cite{Honda_I} for
a different approach to these questions). The main source of
variations lies in boundary conditions since those most
frequently used in \cite{Giroux_2000} are not convenient for our
purposes.

In the following, each time we consider a torus bundle,
the fiber of a point $t$ will be denoted by $T_t$ and $T$ will
denote a torus.

A suspension $\sigma$ on a torus $T$---that is, a 
non-singular foliation admitting a simple closed curve
intersecting all leaves transversely---determines
a line $D(\sigma)$ in $H_1(T,\R)$ called the asymptotic direction
of $\sigma$ and defined in \cite{Schwartzman}: the real cycles
$[L_t]/t$ where $L_t$ is an orbit of length $t$ closed by a
minimizing geodesic segment have a common limit up to sign in
$H_1(T,\R)$ when $t$ goes to infinity and the line $D(\sigma)$
spanned by this limit is independent of the starting point and
Riemannian metric. If $\sigma$ is linear then $D(\sigma)$ is the
direction of $\sigma$ and if $\sigma$ has a periodic orbit then
$D(\sigma)$ is generated by its homology class for any
orientation.  We say that $D(\sigma)$ is rational if it contains
a non-zero point of $H_1(T,\Z)$.

If $\xi$ is a (positive) contact structure on
$T\times I$---where $I$ is an interval---such that each $\xi T_t$
is a suspension then the function
\[
\left(
\begin{array}{ccl}
I &\to& P(H_1(T,\R))\simeq \R/\pi\Z \\
t &\mapsto & D(\xi T_t)
\end{array}
\right)
\]
is non-increasing for the orientation of $H_1(T,\R)$ inherited
from the orientation of $T$ via the intersection form.

\begin{rem}
  \label{rem:tourne} 
If there exists $t$ such that $D(\xi T_t)$ is not rational or 
$\xi T_t$ is conjugated to the suspension of a rotation then
$t\mapsto D(\xi T_t)$ is non-constant. 
\end{rem}

When a $\xi$--convex torus $T$ in a contact manifold has a Legendrian
fibration over the circle then we say that this fibration is a
ruling of $T$ or that $T$ is ruled by this fibration.
In this case $\xi T$ has smooth singularity circles and 
$D(\xi T)$ is defined to be the line in $H_1(T)$ spanned by the
homology class of a singularity circle with any orientation. 

The characteristic foliation of a torus is said to be admissible
if it is linear or ruled with two circles of singularities 
or a suspension divided by two curves.

Tight contact structures on toric annuli, solid tori, lens spaces
and torus bundles with admissible boundary are made of contact
structures printing suspensions on each torus and orbit flips
which we now define (an explicit model is given in
\cite{Giroux_2000}[Section~1.F]).

\begin{defn}
  A contact structure $\xi$ on $T\times [0,1]$ has an orbit
flip with homology class $\pm d$ if all the $\xi T_t$ are divided
by $2k$ circles with total homology class $\pm 2d$ when they are
all given the same orientation, $\xi T_t$ 
is a suspension with $2k$ periodic orbits for
every $t\neq 1/2$, $\xi T_{1/2}$ is ruled 
and the periodic orbits of $\xi T_0$ and $\xi
T_1$ have opposite orientations.
\end{defn}

Note that, in general, the homology class of an orbit flip is
defined up to sign but in our study of contact structures on
Seifert manifolds one of the two possible classes will
intersect the fiber class positively and we will call it
\emph{the} homology class of the orbit flip.

\begin{defn}
  A contact structure on $T\times[0,1]$ is in normal form 
  with flip locus $\{\pm d_1,\dots,\pm d_n\}$
if $\xi T_0$ and $\xi T_1$ are admissible foliations and if 
there exist $t_1<\dots<t_n\in (0,1)$ such that: 
\begin{itemize}
  \item	$\xi_{|T\times (0,t_1)}$ and $\xi_{|T\times (t_n,1)}$
	print a suspension on each $T_t$;
  \item for every $i$, $\xi_{|T\times (t_i,t_{i+1})}$ 
	prints a suspension on each $T_t$ with non-constant asymptotic
	direction;
  \item for every $i$ there exists a neighborhood $J$ of $t_i$ such that 
	$\xi_{|T\times J}$ has an orbit flip with homology class
	$\pm d_i$.
\end{itemize}
\end{defn}

The following definitions are useful to describe the links
between tight contact structures on toric annuli and the geometry
of $H_1(T,\R)$, its integral lattice $H_1(T,\Z)$ and the
intersection form on them.

\begin{figure}[ht]
 \labellist
 \small\hair 2pt
 \pinlabel $D(\sigma_0)$ [bl] at 275 211
 \pinlabel $D(\sigma_1)$ [tl] at 184 25
 \pinlabel $H_1(T,\Z)$ [l] <3pt,2pt> at 12 1
 \pinlabel {$\P$} [l] <3pt,2pt> at 104 1
 \pinlabel {$\partial E\setminus\mathcal{P}$} [l] <3pt,2pt> at 192 1
 \endlabellist
  \begin{center}
	\includegraphics{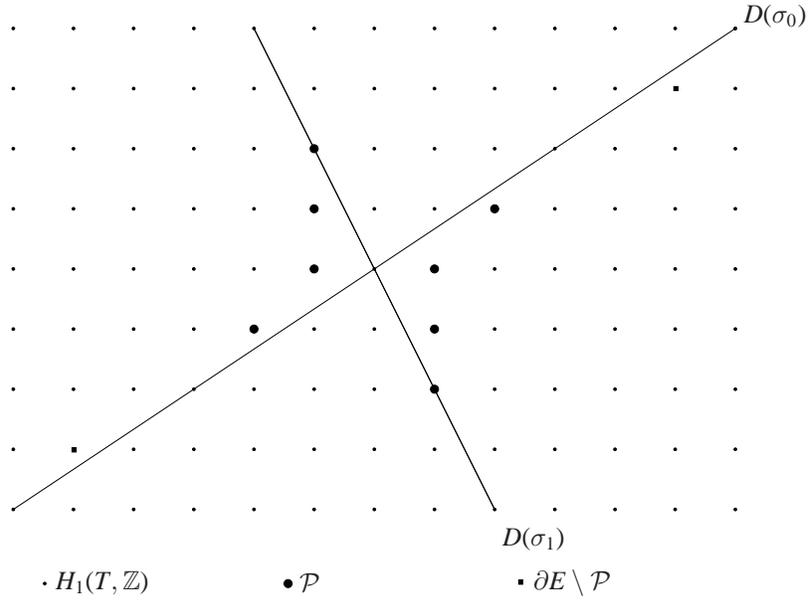}
  \end{center}
  \caption{A Giroux cone and its polygon}
  \label{fig:cone}
\end{figure}

\begin{defn}
 Let $\sigma_0$ and $\sigma_1$ be admissible foliations on $T$.
 The Giroux cone $\mathcal C(\sigma_0,\sigma_1)$ of
 $(\sigma_0,\sigma_1)$ is the cone without vertex bordered
 on the left by $D(\sigma_0)$ and on the right by $D(\sigma_1)$.
 The boundary lines---deprived of 0---are in the cone if and only if the
 corresponding foliations have a dividing set. 
 
 Let $E$ be in each connected component of $\mathcal C$ the
 convex hull of $\mathcal C(\sigma_0,\sigma_1)\cap H_1(T,\Z)$.
 An edge of a part of a lattice is a maximal subset of aligned points.
 The Giroux polygon $\P(\sigma_0,\sigma_1)$ of
 $(\sigma_0,\sigma_1)$ is the set of integral homology classes
 which belong to a finite length edge of 
 $\partial E$---or is the intersection of two infinite edges---, ordered by
 the intersection form from right to left. We denote by 
 $\partial \P$ the (possibly empty) set of its extremal
 points.

 If $\xi$ is a contact structure on $T\times[0,1]$ we also define
 $\mathcal C(\xi)=\mathcal C(\xi T_0,\xi T_1)$ and similarly for
 $\P$.
See figure \ref{fig:cone} for an example where $\sigma_0$ has no
dividing set.
\end{defn}

Theses cones naturally have two connected components and
everything is symmetric with respect to the origin but
in the context of Seifert manifolds we will always consider the
component made of homology classes intersecting
the fiber class positively and call \emph{this component} the Giroux cone.

The following lemma is a special case of
\cite{Giroux_2000}[Lemma~3.34] and can also be proved by
exhibiting an annulus transverse to the $T_t$'s and not
satisfying the Giroux criterion \cite{Giroux_2001}[Theorem~4.5 a)].
We call a element of $H_1(T,\Z)$ simple if it has a simple
curve representative or, equivalently, if it is not a multiple of
an other homology class.

\begin{lem}
\label{lemma:retournement_vrille}
If a contact structure on a toric annulus has an orbit flip
whose homology class is non-simple and not on the boundary of the
corresponding Giroux cone then it is overtwisted.
\end{lem}

For all contact structures printing admissible foliations on the
boundary there is a well defined relative
Euler class once $\sigma_0$ is oriented since such an orientation
fixes a (co)orientation for any (positive) contact structure printing
$\sigma_0$ on $T_0$.
 We denote by $\tight_0(\sigma_0,\sigma_1)$ the set of tight contact structures
with zero torsion printing $\sigma_0$ and $\sigma_1$ on the
boundary components. In the following, a
$\partial$--isotopy is an isotopy relative to the boundary.

\begin{thm}[E Giroux \cite{Giroux_2000}]
  \label{thm:classif_tore_epais}
Let $\sigma_0$ and $\sigma_1$ be admissible foliations on $T$
with $\sigma_0$ oriented and $\P:=\P(\sigma_0,\sigma_1)$.
\begin{enumerate}[a)]
\item Every contact structure in $\tight_0(\sigma_0,\sigma_1)$ is
  $\partial$--isotopic to a contact structure in normal form whose
  flip locus is in $\P$.
\item A contact structure in normal form on $T\times[0,1]$
  which prints $\sigma_0$ on $T_0$ and $\sigma_1$ on $T_1$ is
  tight if and only if its flip locus is contained in 
  $\P$. It is universally tight if and only if its flip locus is
  contained in $\partial\P$.
\item If $D(\sigma_0)\neq D(\sigma_1)$ then
  two contact structures in $\tight_0(\sigma_0,\sigma_1)$ are
  $\partial$--isotopic if and only if they have the same relative
  Euler class.
\item If two contact structures on $T\times[0,1]$ print on each $T_t$ 
	suspensions whose asymptotic directions do not
	cover the whole projective line and if they coincide on the boundary
	then they are $\partial$--isotopic.
\item Any tight contact structure with $\xi$--convex boundary on
	$T\times[0,1]$ is isotopic to a contact structure with
	$T_{1/3}$ (resp. $T_{2/3}$) divided by two curves parallel to
	the dividing set of $T_0$ (resp. $T_1$).
\end{enumerate}
\end{thm}

In the preceding theorem, the correspondence with
\cite{Giroux_2000} is: ~a) and~b) are in Proposition~1.8,
c) is in Theorem~1.5, d) is Theorem~3.3 and~e) follows from
Proposition~3.22.

If a contact structure on $T\times[0,1]$ is in normal form on
$T\times[0,1/2]$ and $T\times[1/2,1]$ with flip loci $R$ and $R'$
then the flexibility lemma and
Theorem~\ref{thm:classif_tore_epais}~d) can be used to prove that
$\xi$ has a normal form with flip locus $(R\cup
R')\setminus(R\cap R')$ so
Theorem~\ref{thm:classif_tore_epais}~b)
can be used to get gluing results. This strategy is used in
Subsection~\ref{SS:collage}.

\subsubsection*{Finite covers}

We now turn to finite covering maps between toric annuli. If
$\rho$ is such a map and $d$ is a homology class represented by a
finite collection $C$ of embedded oriented circles in $T$ then
$\rho^{-1}(C)$ is a finite collection of oriented embedded
circles whose total homology class will be denoted by $\rho^*d$. 

If $\xi$ is a contact structure in normal form on $T\times[0,1]$
with flip locus $\{\pm d_1,\dots,\pm d_n\}$ and $\rho$ is a
finite covering map preserving the product structure of $T\times[0,1]$
then $\rho^*\xi$ has orbit flips with homology classes
$\pm \rho^*d_1,\dots,\pm \rho^*d_n$ so we can use
Lemma~\ref{lemma:retournement_vrille} to detect overtwisted
covers. 

This fact is used
in Section \ref{S:charac} but we can already give a general
corollary: using Theorem~\ref{thm:classif_tore_epais} and
Lemma~\ref{lemma:retournement_vrille} we answer
the question of \cite{Honda_II}[page 97]. Given a virtually
overtwisted contact structure $\xi$ on $T^2\times[0,1]$ with each
boundary component divided by two curves, Honda
asks which covering spaces $\R^2/(m\Z\times n\Z)\times[0,1]$ are
overtwisted. Denote by $\xi(m,n)$ the lifted contact structure.
We now explain how to compute $n_0,m_0$ such that $\xi(m,n)$ is
overtwisted whenever $n\ge n_0$ and $m\ge m_0$. This will leave a
finite number of $\xi(m,n)$ which have to be analyzed directly,
using Theorem~\ref{thm:classif_tore_epais} b).

The crucial point is that $H_1(T^2)$ has a preferred basis and a
euclidean structure in addition to its lattice and intersection
form. According to Theorem~\ref{thm:classif_tore_epais} a) and
b), $\xi$ is isotopic to a contact structure in normal form whose
flip locus is non-empty and contained in
$\mathcal{P}\setminus\partial\mathcal{P}$. Let $d$ be a class in
this flip locus. Let $L_h$ and $L_v$ be affine lines containing
$d$ and directed by the two basis vectors.

Here we assume that $\partial\mathcal C$ has no component
parallel to the axes (the special case we neglect can be dealt with using
the same methods). One can then check that there exist 
$a\in \partial \mathcal C\cap L_h$ and $b\in \partial \mathcal C\cap L_v$
such that the line
$(ab)$ divides $H_1(T^2,\R)$ into two open half-planes, one
containing 0 and the other one containing $d$.
Let $l_h$ (resp. $l_v$) be the distance between $d$ and $a$
(resp. $b$) and set $n_0=\lfloor 1/l_h\rfloor+1$ and 
$m_0=\lfloor 1/l_v\rfloor+1$.
If $n\ge n_0$ and $m\ge m_0$ then the lifted Giroux cone 
contains integer points $a'\in\rho^*[ad]$ and
$b'\in\rho^*[bd]$ proving that $\rho^*d$ is neither on the Giroux
polygon of $\xi(m,n)$ nor on the boundary of its Giroux cone
so that $\xi(m,n)$ is overtwisted according to 
Theorem~\ref{thm:classif_tore_epais}~b) and 
Lemma~\ref{lemma:retournement_vrille}. 

\section{Contact structures on solid tori}
\label{S:toric_quotients}

In this section we explain how the results of the previous
section extend to results on solid tori (here again, everything
comes from \cite{Giroux_2000}) and gather some results about
contact structures with negative twisting number on solid tori
which directly use the classification of tight contact
structures.

\subsection{Classification results}

Let $W$ be an oriented solid torus with a meridian class $M\in
H(\partial W)$, $\sigma$ be an oriented admissible foliation on
$\partial W$ and denote by $\tight(\sigma)$ the set of tight
contact structures on $W$ printing $\sigma$ on the boundary. We
see $W$ as $\S^1\cup(T\times (0,1])$ which is a quotient of
$T\times[0,1]$. A contact structure on $W$ which is transverse to
$\S^1$ can be lifted to $T\times[0,1]$ and is said to be in normal
form if this lift $\xitilde$ is in normal form.
Using this construction $\xitilde$ prints on $T_0$ a linear
foliation $\sigma_M$ whose direction is spanned by $M$ in
$H_1(\partial W)$. The flip locus of $\xi$ is then defined as the
flip locus of $\xitilde$. We denote by $L_M$ the set of integer homology
classes in $\mathcal{C}$ whose intersection with $M$ is $\pm 1$.

\begin{thm}[E Giroux {\cite{Giroux_2000}[Lemma~4.2]}]
  \label{thm:classif_tore_plein}
Let $W$ be a solid torus, $\sigma$ an oriented admissible
foliation on $\partial W$, $\sigma_M$ a linear foliation coming
from meridian disks and $\P$ the Giroux polygon $\mathcal
P(\sigma_M,\sigma)$.
\begin{enumerate}[a)]
\item Every $\xi\in\tight(\sigma)$ is $\partial$--isotopic to a
  contact structure in normal form whose flip locus is in
  $\P$.

\item A contact structure $\xi$ in normal form on $W$ and
  printing $\sigma$ on the boundary is tight if
  and only if its flip locus is contained in 
  $\P\cup L_M$. It is universally tight
  if and only if its flip locus is contained in 
  $\partial\P\cup L_M$.

\item Two contact structures in $\tight(\sigma)$ are
  $\partial$--isotopic if and only if they have the same relative
  Euler class. 
\end{enumerate}
\end{thm}

\begin{cor}
  \label{cor:exist_tore}
  Let $W$ be a solid torus, $\sigma$ an oriented admissible foliation on 
$\partial W$ and $\xi\in\tight(\sigma)$. Denote by $\sigma_M$ a
linear foliation on $\partial W$ spanned by meridian disks.
There exists a
boundary-parallel torus divided by two curves with homology class
$d\in H_1(\partial W,\Z)$ if and only of $d$ is a simple class in
$\mathcal C(\sigma_M,\sigma)$.
\end{cor}

When $\sigma$ is an admissible foliation with dividing curves,
there is only a finite number of $\partial$--isotopy classes of tight
contact structures in $\tight(\sigma)$ which can be computed
using the above theorem.

\begin{cor}
\label{cor:nb_tore_plein}
  With the same notations as above, if $\sigma$ 
is a $\xi$--convex admissible foliation then
$\P$ has a finite number of edges $e_1,\dots,e_k$ and
\[
\Card\bigl(\tight(\sigma)\bigr)=\prod_{i=1}^k \Card(e_i).
\]
\end{cor}

For example, in figure \ref{fig:cone} with $\sigma_0=\sigma_M$ and
$\sigma_1=\sigma$, we have six tight contact structures.
Note that the $\Card(e_i)$ in the preceding corollary can be
computed using the continued fraction expansion of the slope of
$D(\sigma)$  in a suitable basis of $H_1(T,\R)$, see
\cite{Giroux_2000}[Section~1.G].

We now apply this classification to Seifert fibered solid tori:

\begin{lem}
\label{lemma:W_univ_tendue}
Let $W=W(\alpha,\beta)$ be a solid torus with a standard Seifert fibration
(see Subsection~\ref{subsection:notations})
and denote by $M$ a meridian class in $H_1(\partial W)$.
Let $\xi$ be a contact structure on $W$ such that $\xi\partial W$
is ruled by vertical curves with two dividing circles.
Let $d$ be the homology class of these circles oriented such that 
$d\cdot M>0$.

Such a $\xi$ is universally tight if and only if it is 
$\partial$--isotopic to a transverse contact structure 
in the interior of $W$.

If $d\cdot M=1$ then 
there is only one $\partial$--isotopy class 
of tight contact structures which coincide with
$\xi$ on $\partial W$. This class contains universally tight
contact structures tangent to the fibers as well as contact
structures positively and negatively transverse (in
the interior of $W$). 

If $d\cdot M>1$ then there are exactly two
$\partial$--isotopy classes of universally tight contact
structures which coincide with $\xi$ on $\partial W$. They contain
contact structures either positively or negatively transverse
but not both and no tangent contact structure. 
\end{lem}

\begin{proof}
A contact structure in normal form is isotopic through contact
structures in normal form to a transverse contact structure if
and only if its flip locus is empty. This potentially leaves two
isotopy classes of contact structures $\xi'$ depending on the
orientation of the suspensions $\xi' T_t$ near the boundary.
Theorem~\ref{thm:classif_tore_plein} c) and a Euler class
computation give one or two isotopy classes depending on the
intersection number $d\cdot M$ as announced.

Since any tangent contact structure can be
perturbed into a positively or negatively transverse contact
structure in the interior of $W$, it remains only to prove that,
when $d\cdot M=1$, there is a tangent contact structure which
coincides with $\xi$ on $\partial W$.

Let $(S,F)$ be a base of $H_1(\partial W)$ such that $F$ is the
fiber class and $M=\alpha S+\beta F$. Because $d\cdot M=1$, we
can write $d$  as $nS+qF$ with $n\beta - q\alpha= 1$ so that
$W(\alpha,\beta)$ can be seen as
\[
\bigl(D^2 \times \R\bigr)\big/\bigl((z,t)
\sim (e^{\frac{2i\pi n}{\alpha}}z,t+1)\bigr)
\]
because the meridian in the above model has intersection $\alpha$ and $-\beta$
respectively with the fibers induced by the $\R$ factor and the section induced by 
$s \mapsto \bigl(e^{-2i\pi q s},-\beta s\bigr)$.
The map 
\[
\left(
\begin{array}{ccl}
 W(\alpha,\beta) &\to& W(\alpha,1)\simeq ST^*D \\
(z,t) &\mapsto &(z,nt)
\end{array}
\right)
\] 
where $D$ is the base of $W$ is a $n$--fold fibered covering map.  
Seen as a covering map from $W$ to $ST^*D$, it can be chosen to
extend the one over $\partial D$ associated to $\xi$ by (the
oriented version of) \cite{Giroux_2001}[Proposition~3.3]---which
is stated in a slightly generalized form as
Proposition~\ref{prop:rev_giroux_general} below. The pull-back
of the canonical contact structure of $ST^*D$ by this covering
map coincides with $\xi$ on $\partial W$ and is tangent to the
fibers.
\end{proof}

\subsection{A gluing lemma}
\label{SS:collage}

\begin{lem}
  \label{lemma:gluing_pour_calcul_e}
  Let $W(\alpha,\beta)=T^2\times[0,1]/\sim$ be a solid torus with
  a standard Seifert fibration and $\xi$ 
  a contact structure on $W$ such that $T_{1/2}$ and $T_1$ are
  divided by two curves intersecting each fiber only once and the
  fiber class is not in the Giroux cone of $\xi$. If
  $\xi_{|T^2\times[0,1/2]/\sim}$ and $\xi_{|T^2\times[1/2,1]/\sim}$ are
  tight then $\xi$ is tight.
\end{lem}

\begin{figure}[ht]
 \labellist
 \small\hair 2pt
 \pinlabel $F$ [r] at 15 115
 \pinlabel $M$ [bl] at 86 158
 \pinlabel $S$ [tl] <2pt,-2pt> at 29 101
 \pinlabel $d$ [l] <2pt,0pt> at 29 72
 \pinlabel $d\hspace{.25pt}'$ [l] <2pt,0pt> at 29 31
 \endlabellist
 \begin{center}
	\includegraphics{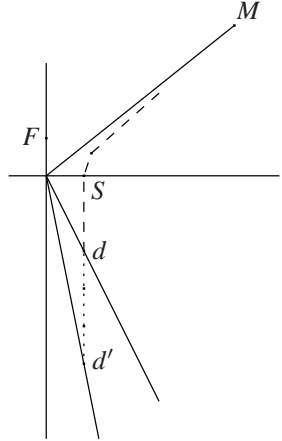}
  \end{center}
  \caption{Giroux polygons for the gluing lemma}
  \label{fig:gluing}
\end{figure}

\begin{proof}
We set
$\xi_1=\xi_{|T^2\times[0,1/2]/\sim}$, $\xi_2=\xi_{|T^2\times[1/2,1]/\sim}$
and we denote by $\P$, $\P_1$ and $\P_2$
the Giroux polygons associated to $\xi$, $\xi_1$ and $\xi_2$.

According to Theorems \ref{thm:classif_tore_plein}~a) and
\ref{thm:classif_tore_epais}~a), we can assume that 
$\xi_1$ and $\xi_2$ are in normal form and that
their flip loci are contained in $\P_1$ and
$\P_2$ respectively. The crucial fact we have to prove
is that $\P =\P_1\cup\P_2$ because of the
assumption on the dividing curves, see figure \ref{fig:gluing}.
This will implies that $\xi$ has a normal form whose flip locus is in
$\P$ and so is tight according to Theorem~\ref{thm:classif_tore_plein} b).

Denote by $F$ and $M$ the fiber and meridian classes oriented
such that $M\cdot F>0$ and denote by $d$ and $\dprime$ the
homology classes of each
dividing curve of $T_{1/2}$ and $T_1$ respectively, oriented such
that their intersection with $F$ is $+1$. We now concentrate
without further notice on the connected component of the 
Giroux cone of $\xi$ which is contained in the half-plane
$\mathcal H = \{h\cdot F>0\}$ and the polygons inside it.
Because $d\cdot F=\dprime\cdot F=1$ and $\dprime\cdot d \geq 0$,
there exists $k\geq 0$ 
such that $d=\dprime+kF$ and, for every $x\in\mathcal H$, the
signed area of the triangle $\dprime x d$ is
\[
(x-\dprime)\cdot(d-\dprime)=(x-\dprime)\cdot kF=
k(x\cdot F-1)\geq 0.
\]
This proves that $\P_2=[\dprime,d]$ and that 
$[\dprime,d]\subset \P$.

In particular, $d$ is in $\P$ and this is exactly what is
needed in order to get $\P=\P_1\cup \P_2$. Indeed, if we denote
by $E$, $E_1$ and $E_2$ the convex hulls of the integer points of
the corresponding cones, what we want is $E=E_1\cup E_2$ so we
need to prove that $E':=E_1\cup E_2$ is already convex. Denote by
$\Delta$ the half-line from $0$ containing $d$.
Since $E_1$ and $E_2$ are convex, we only have to check that, for
every $x\in E_1$ and $y\in E_2$, the point $z:=[x,y]\cap\Delta$
is in $E'$. Since $d$ is in $\P$, $E\cap\Delta=E'\cap\Delta$
and, by convexity of $E$, $z\in E$ so we are done.
\end{proof}

\subsection{Local twisting number calculation}
\label{SS:local_e}

We first recall some arithmetic.
If $x$ is a real number and $a>0$ and $b$ are
relatively prime integers then we say that $b/a$ is a best lower
approximation for $x$ if it is maximum among rational numbers
smaller than $x$ whose denominator is not larger than $a$.
The best lower approximations of $x$ can be read from its
continued fraction expansion:

\[x=[a_0;a_1,\dots]
:=a_0+\frac{1}{a_1+\frac{1}{\ddots}}\]
for a unique sequence of integers $a_i$ such that $a_i>0$ for
$i>0$ and  the last non-zero $a_i$, $i>0$ (if it exists), is greater than 1.

The convergent of order $k$ of $x$ is 
$[a_0;a_1,\dots,a_k]$
and the intermediate fractions of order $k$, $k\ge 2$ of $x$ are the 
$[a_0,a_1,\dots,a_{k-1},a]$ with $1\le a<a_k$.

The best lower approximations of $x$ are exactly the 
convergents and intermediate fractions of even order of $x$ 
with the possible exception of $x$ itself.

\begin{prop}
 \label{prop:local_index}
Let $W$ be a solid torus with a standard Seifert fibration.
Suppose $(S,F)$ is a basis of $H_1(\partial W,\Z)$ such that $F$
is the homology class of the fiber and the meridian circles in
$\partial W$---oriented such that they intersect $F$
positively---have homology class 
$M=\alpha S + \beta F$ with $1\leq\beta<\alpha$ if $\alpha>1$.

If $\xi$ is a tight contact
structure on $W$ such that $\partial W$ is divided by curves 
which can be oriented to have total homology class $d=2nS+2(x-1)F$, $n>0$  
then the following statements are equivalent:

\begin{itemize}
 \item $t(\xi)=-n$;
 \item $(x-1)/n<\beta/\alpha$ and the triangle in $\R^2$ with vertices $(0,0)$,
   $(\alpha,\beta)$, $(n,x-1)$ does not contain any integer
   point whose abscissa is less than $n$ except possibly its
   vertices---in particular $d/2$ is simple;
 \item either $n=1$ and $x-1<\beta/\alpha$ or the following
   conditions hold:
  \begin{itemize}
   \item[(i)]  $x=\lceil n\beta/\alpha \rceil$;
   \item[(ii)]  $n\beta\equiv 1 \mod \alpha$ or $(x-1)/n$ is a best lower
	 approximation of $\beta/\alpha$.
  \end{itemize}
\end{itemize}
\end{prop}

\begin{proof}
Let $\mathcal{C}$ denote the Giroux cone of $\xi$ in
$H_1(\partial W,\R)$.

According to Theorem~\ref{thm:classif_tore_epais}~e),
Corollary~\ref{cor:exist_tore} and the flexibility
lemma,
$t(\xi)=-n$ if and only if $\mathcal{C}$ does not
contain any point $kS+lF$ with $k<n$---in particular $d/2$ is
simple. This proves the equivalence
of the first two points. In addition, the second part of the
second point is clearly always satisfied when $n=1$ so we now
suppose that $n>1$ and we explain the equivalence of the second
and third conditions.

Suppose that the second condition is satisfied.
Set $q=\lceil n\beta/\alpha \rceil-1$ so that $q$ is the 
greatest integer which is smaller than $n\beta/\alpha$.

Recall that, according to Pick's formula (\cite{Pick}, see eg
\cite{the_book} for a proof in English),
any polygon  with integer vertices in $\R^2$ has area
  \[ n_{int}+\frac{1}{2}n_{bd}-1 \] 
where $n_{int}$ is the number of interior integer points and
$n_{bd}$ is the number of boundary integer points.

The point with abscissa $n$ which is the closest to the left
boundary of $\mathcal C$ is $(n,q)$ by definition. Also the
triangle with vertices $(0,0)$, $(n,q)$ and $(n,q-1)$ has area
$n/2>1/2$ so it contains integer points in addition to its
vertices by Pick's formula. Such points necessarily have abscissa
less than $n$ so they cannot be in $\mathcal C$. This implies
that $\mathcal C$ cannot contain $(n,q-1)$, see the left-hand
side of figure \ref{fig:pick}. We conclude that the right
boundary of $\mathcal C$ is generated by $(n,q)$.

\begin{figure}[ht]
 \labellist
 \small\hair 2pt
 \pinlabel $M=(\alpha,\beta)$ [bl] at 115 72
 \pinlabel $(n,q)$ [l] at 86 50
 \pinlabel $(n,q-1)$ [l] at 86 39
 \pinlabel $n$ [tl] <2pt,-2pt> at 86 17
 \pinlabel $M$ [br] at 261 47
 \pinlabel $(n,q)$ [l] at 278 50
 \pinlabel $n$ [tl] <2pt,-2pt> at 278 17
 \endlabellist

 \begin{center}
   \includegraphics{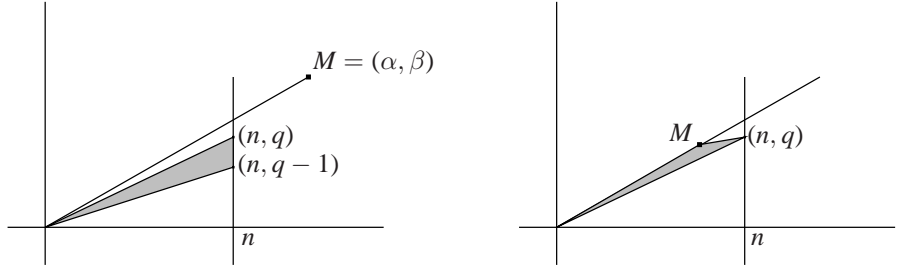}
 \end{center}
  \caption{Computing the multi-index using Pick's formula}
  \label{fig:pick}
\end{figure}

If $n\geq\alpha$ then the triangle with vertices
$(0,0)$, $(\alpha,\beta)$ and $(n,q)$ is inside $\mathcal C$
except for the edge from $(0,0)$ to $(\alpha,\beta)$ which does
not contain any integer point but its vertices, see the
right-hand side of figure \ref{fig:pick}. The area of this triangle is
$(n\beta-q\alpha)/2$ so that according to Pick's formula
we have $n\beta-q\alpha=1$ so $n\beta\equiv 1 \mod \alpha$.

In particular we get that if $\alpha>1$ then $n\neq \alpha$ and 
$\alpha\!\!\not\!|\; n\beta$ so $q=\lfloor n\beta/\alpha\rfloor$ and
$x=\lceil n\beta/\alpha\rceil$.

If $n<\alpha$ then the condition above on integer points in
$\mathcal C$ is exactly the best lower approximation condition.

The fact that the third condition implies the second one is
analogous.
\end{proof}

\section{Characterization of transverse contact structures}
\label{S:debut_charac}

In this section we describe geodesible contact structures on
torus bundles and we prove Theorem~\ref{thm:charac} from the
introduction except for the statement that universally tight
contact structures with negative twisting number on Seifert
manifolds are isotopic to transverse ones which is deferred to
Section~\ref{S:charac}. 

\subsection{Torus bundles}
\label{torus_bundles}

\begin{prop}
\label{prop:charac_t3a}
  Let $A$ be a matrix in $\SL_2(\Z)$ with $\tr(A)>2$.

  A contact structure on $T^3_A$ is isotopic to one which is
transverse to the foliation of Example~\ref{example:t3a}
if and only if it is universally tight and has zero torsion.
There is exactly one isotopy class of such contact structures on
$T^3_A$.

Every geodesible contact structure on $T^3_A$ is symplectically
fillable.
\end{prop}

This classification is a direct consequence of the previously known
classification results contained in \cite{Giroux_2000} (see also
\cite{Honda_I}) which we will recall. The fillability
comes from the classification and a theorem by F~Ding and H~Geiges.

For any function $\theta:\R\to\R$, the Pfaff equation
\[\cos\theta(t)dx_1+\sin\theta(t)dx_2=0\] defines a contact
structure on $T^3_A$ if and only if $\theta'<0$ and
\[A^*(\cos\theta(t+1)dx_1+\sin\theta(t+1)dx_2)\wedge
(\cos\theta(t)dx_1+\sin\theta(t)dx_2)=0.\] This contact structure
is denoted by $\zeta(\theta)$. The following proposition explains
the special role of these $\zeta(\theta)$. It is proved in
\cite{Giroux_2000}[Corollary~3.8] for the first part and the
second part is obtained from the proof of
\cite{Giroux_2000}[Lemma~3.2].

\begin{prop}
  \label{prop:classif_rotative_t3a}
  Let $A$ be a matrix in $\SL_2(\Z)$ with $\tr(A)>2$.
  On $T^3_A$, any contact structure whose lift $\xitilde$ to 
  $T^2\times\R$ prints a suspension on each $T_t$ 
  is isotopic to some $\zeta(\theta)$.

  If $\bigcup_{t\in[0,1]}D(\xitilde T_t)$ is not the whole
  projective line then $\theta(1)-\theta(0)>-2\pi$.
\end{prop}

The theory of normal forms explained in Section
\ref{S:torus_bundles} can also be used to study tight
contact structures on torus bundles over the circle.
In this text we only need to consider the case of universally
tight contact structures on torus bundles with hyperbolic
monodromies. The following theorem comes from
\cite{Giroux_2000}[Theorem~1.3] and its proof in Section~4.D.

\begin{thm}
  \label{thm:classif_t3a}
  Let $A$ be a matrix in $\SL_2(\Z)$ with $\tr(A)>2$. 
  \begin{enumerate}[a)]
	\item   A contact structure on $T^3_A$ is universally tight 
	  if and only if it is isotopic to some $\zeta(\theta)$;
	\item If $\xi$ is conjugated to $\zeta(\theta)$ then
	  $\tor(\xi)=\max\Biggl(\N\cap\left[0,\frac{\theta(0)-\theta(1)}{2\pi}\right)\Biggr)$;
	\item Two universally tight contact structures on $T^3_A$ are isotopic if and only if
	  they have the same torsion.
  \end{enumerate}
\end{thm}

\begin{proof}[Proof of Proposition~\ref{prop:charac_t3a}]
Let $\xi$ be a contact structure transverse to the foliation
$\F_A$ of Example~\ref{example:t3a}.
Because each $T_t$ is foliated by $\F_A$, each
$\xi T_t$ is a suspension and the same is true for the lift 
$\xitilde$ of $\xi$ to $T^2\times\R$.

Also $\bigcup_{t\in[0,1]}D(\xitilde T_t)$ is not the whole
projective line because $\xitilde$ is transverse to a fixed direction.
According to Proposition~\ref{prop:classif_rotative_t3a}, $\xi$ is
isotopic to some $\zeta(\theta)$ with $\theta(1)-\theta(0)>-2\pi$
so it is universally tight and has zero torsion according to
Theorem~\ref{thm:classif_t3a}.

Conversely one can construct a function $\theta$ such that
$\zeta(\theta)$ is transverse to $\F_A$ and all universally tight
contact structures with zero torsion on $T^3_A$ are isotopic to
this $\zeta(\theta)$ according to Theorem~\ref{thm:classif_t3a}.

The symplectic filling is constructed in \cite{Ding_Geiges}[Theorem 1].
\end{proof}

\subsection{Symplectic fillings for transverse contact structures}

Given any Seifert manifold $V$ there is a symplectic manifold
$(W,\omega)$ such that $\partial W=V$ (as oriented manifolds) and
$\ker \omega$ on $V$ is tangent to the fibers. This symplectic manifold
is a (weak convex) filling of $(V,\xi)$ for any transverse
contact structure $\xi$ on $V$.

The existence of $(W,\omega)$ can be deduced from 
\cite{McCarthy_Wolfson}[Theorem 2.1]---as observed by P~Lisca
and G~Mati\'c in \cite{Lisca_Matic}---or, more elementarily, from
the main theorem of \cite{Klaus_Frederica}. Indeed, if
$(B,\omega_B)$ is a symplectic orbifold with isolated cyclic
singularities and $V\to B$ is a circle bundle over $B$ (in the
sense of Thurston \cite{Thurston_chap_13}) then one can consider
the associated disk bundle $D\to B$ where each fiber is equipped with (a
quotient of) the symplectic form $d(\frac{1}{2}r^2d\theta)$.
Following the construction explained eg in \cite{MDS}[Theorem
6.3] we get a symplectic orbifold which can be resolved using
\cite{Klaus_Frederica} to get the desired filling. If $B$ is a
2--orbifold like in our case then the resolution is completely
explicit.

\subsection{Universal tightness}
\begin{prop}
\label{prop:geod_univ_tendue}  
If $(V,\xi)$ is a Seifert manifold with a transverse contact
structure then its universal cover $(\Vtilde,\xitilde)$ is
$\R^3$ or $\S^3$ with its standard tight contact structure.
\end{prop}

\begin{proof}
  The universal cover $\Vtilde$
of $V$ is either $\R^3=\R^2\times\R$ or $\\S^2\times\R$ or $\\S^3$
where the Seifert structure lifts to the foliation by $\R$ in the
first two cases and $\\S^3$ can have any of its Seifert structures
(see Example~\ref{ex:s3}).

Here $\Vtilde$ cannot be $\S^2\times\R$ since $V$ would
be covered by $\S^2\times \S^1$ which has no transverse contact
structure. This fact is contained in Theorem~\ref{thm:existence}
but there is a direct argument. Indeed, $\S^2\times \S^1$ is
the quotient of $\R/\Z\times\R/\Z\times [0,1]$ by 
$(x,y,t)\sim(x',y,t)$ when $t=0$ or $1$ which is foliated by
$\{\cdot\}\times\R/\Z\times\{\cdot\}$. If $\xi$ is a transverse contact
structure then it can be lifted to a contact structure
$\xitilde$ which prints a suspension on each torus $T_t$. The
foliations $\xitilde T_0$ and $\xitilde T_1$ are conjugated to
the suspension of the same rotation so the asymptotic direction
of $\xitilde T_t$ in $H_1(T^2,\R)=\R\times\R$ is
not constant according to remark \ref{rem:tourne} and goes from
$\R\times\{0\}$ to itself. This contradicts the fact that it
never is $\{0\}\times\R$ which is the fiber direction.

If $\Vtilde$ is $\R^3$ then $(V,\xi)$ is covered by $\R^2\times \S^1$
equipped with a contact structure transverse to the $\S^1$ factor.
So we can decompose $\R^3$ as $\R^2\times\R$ so that the lifted
contact structure $\xitilde$ is transverse to the $\R$ factor and invariant
under integral translations in this direction. According to 
\cite{Giroux_2001}[Section 2.B.c], $\xitilde$ is tight
(this is proved using explicit contact embedding of any ball inside 
the standard $\R^3$ where the Bennequin theorem \cite{Bennequin} 
is available).
According to a theorem of Y~Eliashberg \cite{Eliashberg_20_ans},
$\xitilde$ is the standard contact structure on $\R^3$.

Suppose now that $\Vtilde$ is $\S^3$. We have to prove that 
any transverse contact structure $\xi$ on $\S^3$ is
isotopic to the standard contact structure on $\S^3$, which is
tight according to the Bennequin theorem.

Let $p:T^2\times[0,1]\to \S^3$ be the projection introduced in
Example~\ref{ex:s3} and $\xihat$ the lifted contact structure on
$T^2\times[0,1]$. Up to an isotopy of $\xi$ among transverse
contact structures, we can assume that $\xi$ coincides with the
standard contact structure $\xi_0$ of $\S^3$ in a neighborhood of
$L_0\cup L_1$. According to
Theorem~\ref{thm:classif_tore_epais}~d) applied 
to $\xihat$ and $p^*\xi_0$, $\xi$ is isotopic to the
standard contact structure of $\S^3$.
\end{proof}

\subsection{Twisting number}

\begin{prop}
  \label{prop:geod_negative}
  Transverse contact structures on Seifert manifolds have
negative twisting numbers.
\end{prop}

\begin{proof}
Let $V$  be a Seifert manifold and $\xi$ a transverse contact structure on
$V$.
  We first remark that a contact structure which is covered by a
contact structure with negative twisting number has
negative twisting number. Indeed, suppose $\Vhat \to V$ is a
$k$--fold fibered covering map, $\xihat$ is lifted from $\xi$,
$f$ is a regular fiber in $V$ and $\skew{3.5}\hat f$ the lifted
regular fiber of $\Vhat$. For every isotopy bringing $f$
to a $\xi$--Legendrian curve $L$ in $V$ there is a lifted
isotopy which brings $\skew{3.5}\hat f$ to a 
$\xihat$--Legendrian curve $\Lhat$ in $\Vhat$ and
$t(\Lhat)=k t(L)$ so if $L$ has non-negative twisting number then
so has $\Lhat$.
Note that in general $t(\xihat)$ can nonetheless be higher than 
$t(\xi)$ because of Legendrian curves which are not lifted from
curves in $V$ so Theorem~\ref{thm:e} cannot be deduced from the
circle bundle case and indeed exhibits a much richer behavior in
the general case.

If the base orbifold $B$ of $V$ 
is covered by a smooth surface then $V$ is covered by a circle
bundle $\Vhat$ and we conclude using the preliminary remark and
\cite{Giroux_2001}[Proposition 2.4c].
 If $B$ is not covered by a smooth surface then $V$ is a lens
space whose universal cover is $\S^3$. If this universal cover is
a circle bundle then we can use Giroux's result again but in
general $\S^3$ will have one or two exceptional fibers, see
Example~\ref{ex:s3}. According to Proposition~\ref{prop:geod_univ_tendue} the lifted contact structure
coincides with the standard contact structure on $\S^3$. Regular
fibers of $\S^3$ are then positive torus knots and we can conclude
using Bennequin's inequality. Indeed, if the regular fibers are
$(p,q)$ torus knots then the Seifert framing and the fibration
framing differ by $pq$ (this is the linking number of two
disjoint $(p,q)$ torus knots contained in the same torus) and the
Seifert genus of a $(p,q)$ torus knot is easily seen to be at
most $(p-1)(q-1)/2$, see eg \cite{Rolfsen}[chapter 5] (this is
actually the exact Seifert genus).
Bennequin's inequality applied to any vertical Legendrian curve
$L$ then gives
$t(L)=tb(L)-pq\leq 2g-1-pq \leq -p-q<0.$
\end{proof}

\begin{rem}
\label{remark:framing}
If $V$ is a Seifert manifold such that there exists an isotopy
relative to a regular fiber which does not preserve the canonical
framing of this fiber then every contact structure has infinite
twisting number. Because of the above proposition we are
not interested in those Seifert manifolds in this paper so we can
safely forget about the isotopy when we consider vertical curves
and still have a canonical framing.
\end{rem}

In some papers, the twisting number is defined to be zero
whenever it is non-negative. This discrepancy has no impact in
the present paper since we will study almost exclusively negative
twisting numbers and we stick to the definition of
E~Giroux in \cite{Giroux_2001} because of the following lemma.

\begin{lem}
\label{lemma:vrille_e}
  Overtwisted contact structures on Seifert manifolds have
infinite twisting number.
\end{lem}

This is \cite{Giroux_2001}[Proposition 2.4.b]
which was written in the context of circle bundles but the
statement and proof are the same (take the connected sum of a
vertical Legendrian curve and a Legendrian unknot with positive
Thurston--Bennequin invariant near an overtwisted disk to increase
twisting number arbitrarily).

This lemma tells us that in order to prove that a contact
structure $\xi$ is tight, it is sufficient to prove that
$t(\xi)<0$. Also the original proof of Bennequin's theorem that the
standard contact structure on $\S^3$ is tight consists in proving that
any Legendrian unknot has non-positive Thurston--Bennequin
invariant and this is completely equivalent to the fact that this
contact structure has negative twisting number for the circle
bundle structure of $\S^3$. However
there exist tight (and even universally tight) contact
structures on circle bundles with zero twisting number,
see \cite{Giroux_2001}.

\section{Existence criterion on Seifert manifolds}
\label{S:existence}

In this section we prove Theorem~\ref{thm:existence} and its
corollary from the introduction.

\subsection{Arithmetic criteria}

We first explain how the corollary stated in the introduction as well
as the following one are deduced from Theorem~\ref{thm:existence}.

\begin{cor}
  Let $V$ be a Seifert manifold. The following statements are
equivalent:
\begin{itemize}
  \item $V$ has a transverse contact structure $\xi$ with
	$t(\xi)=-1$;
  \item $V$ has a contact structure $\xi$ with
	$t(\xi)=-1$;
  \item $e_0(V)\leq 2g-2$.
\end{itemize}
\end{cor}

In order to get concise statements, we recall a definition coming from \cite{JN_2} and
\cite{Naimi}:
  a tuple $(\gamma_1,\dots,\gamma_r)\in\left(\Q\cap(0,1)\right)^r$ 
  is \emph{realizable} if $r\geq 3$ and if there
exist relatively prime integers $0<a<m$ such that---possibly
after reordering the tuple---we have
$\gamma_1<a/m$, $\gamma_2<(m-a)/m$ and $\gamma_i<1/m$ for every 
$i\geq 3$.
If $V$ is a Seifert manifold we set
$\Gamma(V)=(1-\beta_1/\alpha_1,\dots,1-\beta_r/\alpha_r)$.

Using the equivalence of (i) and (iii) in Theorem~\ref{thm:existence} we only have to prove the following 
purely arithmetic fact:

\begin{prop}
  Let $g$ and $r$ be non-negative integers, $b$ an integer and
  \[(\alpha_1,\beta_1),\dots,(\alpha_r,\beta_r)\] pairs of integers
  such that $1\leq \beta_i <\alpha_i$ for every $i$. We set
  $e_0=-b-r$ and $e=-b-\sum\beta_i/\alpha_i$. The following statements are equivalent:
  \begin{itemize}
	\item[(a)] There exist a positive integer $n$ and integers
	  $x_0,\dots,x_r$ such that 
	  \[(x_0-1)/n<b,\;\; (x_i-1)/n<\beta_i/\alpha_i \mbox{ for
	  } i>0 \mbox{ and } \sum x_i=2-2g;\]
	\item[(b)] One of the following holds:
	  \begin{itemize}
		\item[(i)] $e_0\leq 2g-2$;
	    \item[(ii)] $g=0$, $r\leq 2$ and $e<0$;
		\item[(iii)] $g=0$, $r\geq 3$, $e_0=-1$ and 
		  $(1-\beta_1/\alpha_1,\dots,1-\beta_r/\alpha_r)$ is
		  realizable.
	  \end{itemize}
  \end{itemize}
  In addition, (a) is true with $n=1$ if and only if 
  $e_0\leq 2g-2$.
\end{prop}

\begin{proof}
Remark first that if (a) is true then $x_0\leq nb$
and $x_i\leq n$ for every $i\geq 1$. Also $2-2g=\sum x_i$
so $2-2g\leq n(b+r)=-ne_0$ so that $ne_0\leq 2g-2$.

\paragraph{Non-zero genus case}

Suppose $g>0$. If (a) is true then the above remark gives 
$e_0\leq (2g-2)/n\leq 2g-2$.

Conversely if (b) is true then $e_0\leq 2g-2$ and we can choose
$n=1$, $x_0\leq b$ and $x_1,\dots x_r \leq 1$ such that 
$\sum x_i=2-2g$.

\paragraph{Zero genus case}

We now consider the case $g=0$. The same argument as above shows
(a) is true with $n=1$ if and only if $e_0\leq 2g-2=-2$.

Suppose that (a) is true with $n>1$.
The preliminary remark gives us $e_0\leq (2g-2)/n$ which is
negative so $e_0\leq -1$ because $e_0$ is an integer.

If $e_0\leq -2$ then we conclude using the first remark. Suppose
now that $e_0=-1$, so that $b=1-r$.

If $r=0$ then $b=1$ so $e=-b<0$.
If $r\in\{1,2\}$
 \[2=\sum_{i=0}^r x_i < nb +\sum_{i=1}^r\left(\frac{n\beta_i}{\alpha_i}+1
\right)=-ne+r\]
so $e<(r-2)/n\le 0$ and $e<0$. 

Conversely if $e<0$ then $b\geq 1-r$ and if $r=0$ 
then we can choose
$n=2$, $x_0=2$, if $r=1$ then $e<0$ implies the existence of $n$
such that $e<-1/n$ and we choose $x_0=nb$ and $x_1=2-x_0$.
If $r=2$ then $e<0$ means $-b-\beta_2/\alpha_2<\beta_1/\alpha_1$
so there exist positive $k$ and $n$ such that 
$-b-\beta_2/\alpha_2<k/n<\beta_1/\alpha_1$ and we choose
$x_0=nb$, $x_1=k+1$ and $x_2=2-x_0-x_1$.

The only remaining case is $g=0$, $e_0=-1$ and $r\geq3$.
If (b) is true then 
we choose $n=m$, $x_0=m(1-r)$, $x_1=m-a+1$, $x_2=a+1$ and
$x_i=m$ for every $i\geq 3$.

Conversely if (a) is true then $x_0\leq nb$ and we can replace
$x_0$ by $nb$ and $x_1$ by $x_1-(nb-x_0)$ without losing
anything so we can assume that $x_0=nb$. We claim that $x_i\geq
2$ for every positive $i$. Indeed if this is not true then
\[2=\sum_{i\geq 0}x_i<nb+2+n(r-1)=2\]
which is absurd. So $x_i\geq 2$ and we conclude using the following  lemma.
\end{proof}

\begin{lem}
  If $r\geq 3$, $n\geq 2$ and $2\leq x_1\leq\dots\leq x_r\leq n$ are
integers such that $\sum x_i = 2+n(r-1)$ then there exist
relatively prime integers $0<a<m$ such that
\[\frac{x_1-1}{n} \geq \frac{m-a}{m},\qquad 
\frac{x_2-1}{n} \geq \frac{a}{m},\qquad \forall i\geq 3,\;
\frac{x_i-1}{n} \geq \frac{m-1}{m}.\]
\end{lem}

\begin{proof}
  First remark that we do not have to care about $a$ and $m$
being relatively prime because we can always divide them by their
greatest common divisor while retaining their relations to the
$x_i$.

Set $k=x_1+x_2-2-n$. The integer $k$ is non-negative because 
\[x_1+x_2=2+n(r-1)-\sum_{i=3}^r x_i \geq
2+n(r-1)-n(r-2)=2+n\]
(we replaced the $x_i$ by $n$ in the sum).

Set $l=n-x_3$. By assumption $l$ is non-negative and it is
not larger than $k$ since
\[2+n(r-1)=\sum x_i=2+n+k+n-l+\sum_{i\geq 4} x_i \leq
2+n+k+n-l+n(r-3).\]

Since the $x_i$ form a non-decreasing sequence, we only need to
show that there exist $a$ and $m$ meeting the conditions related
to $x_1$ and $x_2$ and such that 
$(n-l-1)/n\geq (m-1)/m$. The latter condition is equivalent to
$m\leq n/(l+1)$. 
Put $m=\lfloor n/(l+1)\rfloor$. We know $m\geq 2$ because
$x_1\leq n$ so $x_2\geq 2+k$ and $n-l=x_3\geq x_2$ so $n-l\geq 2+k \geq 2+l$
so $n/(l+1)\geq 2$. We only need to check the existence of
$0<a<m$ meeting the conditions related to $x_1$ and $x_2$. 

In the Euclidean plane $\R^2$ we consider points
$A=(n,x_1-1)$, $B=(n,-(x_2-1))$, $H=(n,0)$ and $H'=(m,0)$
(see figure \ref{fig:thales}).

Denote by $A'$ (resp. $B'$) the intersection point between the
line $x=m$ and the line $(OA)$ (resp. $(OB)$). 
According to the intercept theorem the segment
$[A'B']$ has length $(n+k)m/n\geq m$ therefore it contains a
segment $[A''B'']$ with length $m$ whose extremities have integer
coordinates. Let $m-a$ be the ordinate of $A''$ so that the
ordinate of $B''$ is $-a$. The integer $a$ is such that 
$(x_1-1)/n\geq (m-a)/m$ and $(x_2-1)/n\geq a/m$ and we have
$a>0$ because if $a=0$ then $x_1\geq n+1$, which is absurd.
\begin{figure}[ht]
 \labellist
 \small\hair 2pt
 \pinlabel $A$ [bl] at 111 129
 \pinlabel $H$ [tl] <1pt,-1pt> at 111 86
 \pinlabel $B$ [tl] at 111 16
 \pinlabel $A'$ [br] at 55 105
 \pinlabel $A''$ [l] at 55 101
 \pinlabel $H'$ [tl] <1pt,-1pt> at 55 86
 \pinlabel $B''$ [bl] at 55 58
 \pinlabel $B'$ [tr] at 55 56
 \endlabellist

  \begin{center}
	\includegraphics{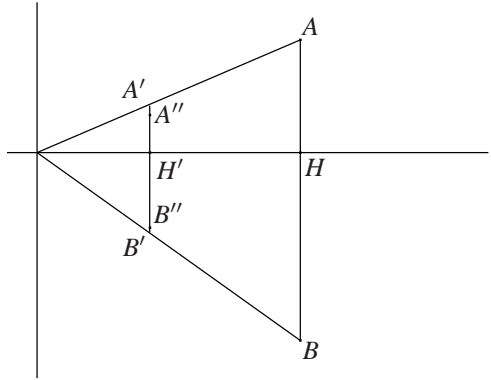}
  \end{center}
  \caption{Intercept theorem and realizability}
  \label{fig:thales}
\end{figure}
\end{proof}

Using the Eliashberg--Thurston perturbation theorem in \cite{ET}
we can recover the following result about foliations:

\begin{cor}[\cite{EHN,JN_1,JN_2,Naimi}]
  A Seifert manifold $V$ with a $C^2$ transverse foliation
satisfies one of the following conditions:
\begin{itemize}
 \item $e_0(V) \leq 2g-2$ and $e_0(-V) \leq 2g-2$;
 \item $g=0$ and $e(V)=e(-V)=0$;
 \item $g=0$, $e_0(V)=-1$ and $\Gamma(V)$ is realizable;
 \item $g=0$, $e_0(-V)=-1$ and $\Gamma(-V)$ is realizable.
\end{itemize}
\end{cor}

\subsection{Normal forms}

Let $R$ be a compact orientable surface of genus $g$ with $r+1$
boundary components and $J$ be a complex structure on $R$
defining its orientation.
The Seifert manifold $V$ with invariants 
\[(g,b,(\alpha_1,\beta_1),\dots,(\alpha_r,\beta_r))\]
is the union of $R\times \S^1$ and $r+1$ solid tori 
$W_0,\dots, W_r$, see Section \ref{subsection:notations}.

Recall from the introduction that a contact structure $\xi$ on $V$ is
in normal form if its restriction to $R\times \S^1$ has an equation
$\cos(n\theta)\lambda+\sin(n\theta)\lambda\circ J$ where
$\lambda$ is some non-singular 1--form on $R$ and $n=-t(\xi)$.
Such a $\xi$ is denoted by $\xi(\lambda,n,\xi_0,\dots,\xi_r)$
where $\xi_i=\xi_{|W_i}$.
Also, the multi-index of this contact structure  is
the collection of the indices of $\lambda$ along the boundary
components of $R$.

\begin{prop}
  \label{prop:exist_normal_form}
  Every contact structure having negative twisting number
on a Seifert manifold 
$V=(R\times \S^1)\cup W_0\cup\dots\cup W_r$ is isotopic to a contact
structure in normal form.
\end{prop}

This proposition follows from ideas in \cite{Giroux_2001} which
we briefly recall. The following lemma is a consequence of the
flexibility lemma:

\begin{lem}
\label{lemma:exist_ruling}
Let $\xi$ be a contact structure having negative twisting
number on a Seifert manifold.
Any $\xi$--convex vertical torus which contains a maximally twisting
vertical Legendrian curve is isotopic through $\xi$--convex surfaces to a
torus ruled by maximally twisting vertical Legendrian curves. Any
$\xi$--convex annulus whose boundary consists of two vertical
Legendrian curves with maximal twisting number is divided by curves
going from one boundary component to the other one and is
$\partial$--isotopic through $\xi$--convex surfaces 
to an annulus ruled by maximally twisting 
vertical Legendrian curves.
\end{lem}

Using this lemma as in \cite{Giroux_2001}[lemma 2.8] one 
proves that any contact structure having negative twisting
number $-n$ is isotopic to a contact structure such that every
fiber in $R\times \S^1$ is Legendrian with twisting number $-n$.
These contact structures over $R\times \S^1$ are all pulled back
from $\xi_R$ on $ST^*R$ by fibered covering maps as explained by
(the oriented version of) \cite{Giroux_2001}[Proposition 3.3]
(which is stated in a slightly generalized form as
Proposition~\ref{prop:rev_giroux_general} below). These covering
maps can be deformed to get a contact structure in normal form.

\subsection{Criterion}

\begin{proof}[Proof of Theorem~\ref{thm:existence}]
  According to Proposition~\ref{prop:geod_negative}, transverse contact
structures have negative twisting number so $(i)\Rightarrow (ii)$. 
  
  We now prove that $(ii)\Rightarrow (iii)$. Let
$\xi=\xi(\lambda,n,\xi_0,\dots,\xi_r)$ be a contact
structure in normal form in the isotopy class of $\xi$ given by
$(ii)$ and Proposition~\ref{prop:exist_normal_form} and let
$(x_0,\dots,x_r)$ be its multi-index. Since normal forms are
completely explicit, we easily see that each $\xi\partial W_i$ is
divided by curves whose total homology class (when they are all
oriented to intersect positively $F_i$) is $2nS_i+2(x_i-1)F_i$.

According to the Poincar\'e--Hopf theorem applied to $\lambda$,
$\sum x_i=2-2g$, and according to Proposition~\ref{prop:local_index}, 
$(x_i-1)/n<\beta_i/\alpha_i$.

 We now prove that $(iii)\Rightarrow (i)$. 
Since $\sum x_i=2-2g$, there exists a non-singular 1--form
$\lambda$ on $R$
whose indices along the boundary components $\partial D_i$ of $R$
are the $x_i$'s. The corresponding contact structure
$\cos(n\theta)\lambda+\sin(n\theta)\lambda\circ J$ on 
$R\times \S^1$  can be extended by
contact structures $\xi_i$ inside each $W_i$ which are positively
transverse, first in a neighborhood of $\partial W_i$ using the
flexibility lemma and then explicitly by a contact structure
printing suspensions on concentric tori around the central fiber
of $W_i$ because $(x_i-1)/n<\beta_i/\alpha_i$ so the asymptotic
directions of these suspensions can go from the half-line spanned
by $M_i$ to the one spanned by $nS_i+(x_i-1)F_i$ without crossing
the fiber direction.
 We then push slightly the
resulting contact structure on $V$ using a Legendrian vector
field which is orthogonal to the fibers over $R$ and zero inside
the $W_i$'s and then use a small isotopy near $\partial R\times
\S^1$ to get a transverse contact structure $\xi$.
 Of course $t(\xi)\geq -n$ and $t(\xi) <0$ according to
Proposition~\ref{prop:geod_negative}.
\end{proof}

\section{Universally tight normal forms}
\label{S:charac}

In this section we end the proof of Theorem~\ref{thm:charac} from
the introduction with the following proposition:

\begin{prop}
\label{prop:univ_tendue_positive}
  Let $V=(R\times \S^1)\cup W_0\cup\dots\cup W_r$ be a Seifert
manifold.
A contact structure $\xi=\xi(\lambda,n,\xi_0,\dots,\xi_r)$ in
normal form on $V$ is universally tight if and only if it can be
cooriented such that each $\xi_i$ is $\partial$--isotopic to a
positively transverse contact structure.
\end{prop}

\begin{proof}
We have already explained in the proof of
Theorem~\ref{thm:existence} that if all the $\xi_i$'s are
positively transverse then $\xi$ is isotopic to a (positively)
transverse contact structure hence universally tight according
to Proposition~\ref{prop:geod_univ_tendue}.

Conversely, suppose that $\xi$ is  universally tight. 
According to Proposition~\ref{prop:local_index}, all the
$\xi_i\partial W_i$ are divided by two circles whose 
homology class---when they are oriented to intersect the fibers
positively---will be denoted by $d_i$.

According to Lemma~\ref{lemma:W_univ_tendue}, there are two
things to check. The first step is to prove that all the $\xi_i$'s
are universally tight and this implies that they are all
$\partial$--isotopic to transverse contact structures.
Then if $d_i\cdot M_i=1$ for every $i$ there is nothing more to prove
since all the $\xi_i$'s are isotopic to positively transverse
contact structures. Else if there is some $i_0$ such that
$d_{i_0}\cdot M_{i_0}>1$ then we can coorient $\xi$ such that
$\xi_{i_0}$ is positively transverse and the second step is to prove
that for every $j$ such that $\xi_j$ is negatively transverse we
have $d_j\cdot M_j=1$.

In the following we will use normal forms of contact structures
on toric annuli and solid tori from Sections
\ref{S:torus_bundles} and \ref{S:toric_quotients}. There
should not be any confusion with normal forms defined in the
introduction since they do not live on the same
manifolds. For each $i$ we put $\xi_i$ in normal
form using Theorem~\ref{thm:classif_tore_plein}~a) and then,
according to Theorem~\ref{thm:classif_tore_plein}~b) $\xi_i$ is
universally tight if and only if every $d$ in its flip locus
satisfies either $d\cdot M_i=1$ or $d=d_i$.

Let $(\Vtilde,\xitilde)$ be the universal cover of $(V,\xi)$
and $\rho$ the covering map. According to
Theorem~\ref{thm:existence}, $\Vtilde$ is $\S^3$ or $\R^3$---this
can be seen using $(ii)\Rightarrow(iii)$ and the fact that a
Seifert manifold is covered by $\S^2\times\R$ if and only if its
base is spherical and $e=0$ or using $(ii)\Rightarrow(i)$ and
Proposition~\ref{prop:geod_univ_tendue}.
We will explain in detail the case where $\Vtilde=\S^3$ with no
exceptional fiber---which the subtlest---and indicate briefly
how to deal with the other cases.

{\bf Preliminary observations }
If a Seifert manifold $Y$ is the union of tubular neighborhoods of
two fibers $K$ and $K'$ then it can be seen as a quotient of
a toric annulus $T\times[0,1]$---where $T_0$ projects to $K$ and
$T_1$ to $K'$---and if $Y$ has some contact structure $\xi$ 
with negative twisting number then there are meridian classes 
$M$ and $M'$ coming from $K$ and $K'$ in $H_1(T)$ such that 
$M'\cdot M$, $M'\cdot F$ and $M\cdot F$ are positive. If $\xi$ is transverse to
$K$ and $K'$ then it can be lifted to a contact structure
$\xihat$ on $T^2\times[0,1]$. If $Y$ is $\S^3$ then $M'\cdot M=1$ and
the Giroux polygon of $\xihat$ is $\{M'+M\}$ so every $d$ in the
flip locus of a normal form of $\xihat$ satisfies $M'\cdot d=1$ or
$d\cdot M=1$.

{\bf First case }
Suppose first that $\Vtilde=\S^3$ with no exceptional fiber. 

For the first step we fix some $i$ and we consider a connected
component $\Wtilde$ of $\rho^{-1}(W_i)$. The
complement $\Wtilde'$ of $\Wtilde$ in $\S^3$ is also a solid
torus and we are in the situation of the preliminary
observations.

For every $d$ in the flip locus of $\xi_i$, $\xihat$ has an orbit
flip with homology class $\rho^*d$ so $\rho^*d$ is simple
according to Lemma~\ref{lemma:retournement_vrille} and $\xitilde|_{\Wtilde}$
is in normal form and then we put $\xihat$ in
normal form using Theorem~\ref{thm:classif_tore_epais}~a).

According to Proposition~\ref{prop:local_index}, 
the triangle with vertices $0,M_i,d_i$ doesn't
contain any integer point $d\not\in\{0,M_i\}$  with $d\cdot F<d_i\cdot F$.
 In particular such a $d$ cannot be in the triangle $\Delta$ with vertices
$0, F, d_i$ because the triangle with vertices $0, F, d$ would
then be included in $\Delta$ so it would have area $d\cdot F/2$
smaller than the area $d_i\cdot F/2$ of $\Delta$.
So for every $d\neq d_i$ in the flip locus of $\xi_i$,
$d\not\in\Delta$ and  $\rho^*d\not\in\rho^*\Delta$---because
$\rho^*$ is linear. The point is that the latter triangle contains all
the integer points $a$ such that $M'\cdot a=1$, $\rho^*d_i\cdot a>0$ and
$a\cdot F>0$ so $M'\cdot \rho^*d>1$ and $\rho^*d\cdot M=1$. This implies that
$d\cdot M_i=1$ and $\xi_i$ is universally tight.

For the second step, let $\Wtilde$ be a connected component of
$\rho^{-1}(W_{i_0})$ and $N$ be a toric annulus around $\Wtilde$
containing one connected component of $\rho^{-1}(W_j)$ and no
other connected component of a $\rho^{-1}(W_i)$. The torus
$\partial \Wtilde$ is divided by curves with total homology
class $2\rho^*d_{i_0}$ and the other component of $\partial N$ is
divided by curves with total homology class $2\dprime'$ such that
$\dprime'-\rho^*d_{i_0}=kF$ because these two tori are ruled by
Legendrian fibers with the same twisting number and $k\le 0$
because $\dprime'\cdot \rho^*d_{i_0}\geq 0$. Because of the additivity
property of the indices of $\lambda$ along curves, $d_j\cdot M_j=1$ if
and only if $k=0$.

We isotop
$\xitilde$ in $N$ such that it is negatively transverse 
in the interior of $N$. The contact structure
$\xitilde_{|\Wtilde\cup N}$ is then in normal form with flip
locus $\{\rho^*d_{i_0}\}$. The complement of $\Wtilde\cup N$ in
$\S^3$ is a solid torus $\Wtilde''$ and we are again in the
situation of the preliminary observations.
We put $\xitilde_{|\Wtilde''}$ in normal form and then
$\xihat$ is in normal form and---because
$\rho^*d_{i_0}\cdot M>1$---we get $M''\cdot \rho^*d_{i_0}=1$. So we have 
$0<M''\cdot \dprime'=M''\cdot \rho^*d_{i_0}+kM''\cdot F=1+kM''\cdot F$ hence $-1<k\le 0$ so
$k=0$.

{\bf Second case }
If $\Vtilde$ is $\S^3$ with at least one exceptional fiber then
$V$ is a Lens space.  We see $V$ as $W\cup N\cup W'$ where all
pieces are fibered, $W$ and $W'$ contain at most one exceptional
fiber each, $W_0\subset N$ and $N$ intersects no other $W_i$. 
We are in the situation of the preliminary observations and we
denote by $2d$ and $2\dprime$ the total homology classes of the
dividing sets of $\partial N$ seen in
$H_1(T)$ and oriented such that $d\cdot F>0$ and $\dprime\cdot F>0$.

For the first step we can use the same argument as above except
when $i=0$.  If $n>1$ then according to Proposition~\ref{prop:local_index} and Lemma~\ref{lemma:W_univ_tendue}
$\xi_0$ is universally tight (and isotopic to a tangent contact
structure). If $n=1$ then the Giroux polygon of $\xi_{|N}$ 
is $[d_1,d_2]$ and is included in the Giroux polygon of
$\xihat$. But, still because of the Giroux polygon associated to
$\S^3$ in the preliminary observations, $\xihat$ is universally
tight so $\xi_{|N}$ cannot have an orbit flip in the interior
of its polygon so it is universally tight and so is $\xi_0$
because any cover of $W_0$ is contained in a cover of $N$.

For the second step we can directly use that every $\dhat$ in
the flip locus of $\xihat$ satisfies $M'\cdot \dhat=1$ or 
$\dhat\cdot M=1$ and that $\xi_0$ is isotopic to a tangent contact
structure if and only if $d=\dprime$. 

{\bf Third case }
If $\Vtilde=\R^3$ then all the $W_i$'s are covered
by some $D^2\times\R\subset\R^3$ so all the $\xi_i$'s are universally
tight. 

The second step is analogous to the first case but we only go to
a cover by $\R^2\times \S^1$ and use that all fibered tori in it
have universally tight lifted contact structures. Using that the
analogous of $\xitilde_{|\Wtilde\cup N}$ is universally tight
we directly get that $\rho^*d_{i_0}=\dprime'$.

Note that in the $\R^3$ case we don't need 
$\xi(\lambda,n,\xi_0,\dots,\xi_r)$ to have twisting number $-n$
to get the result.
\end{proof}

\section{Maximal twisting number calculations}
\label{S:calcul_e}

In this section we prove Theorem~\ref{thm:e} from the
introduction. We also note the following corollary of
Theorem~\ref{thm:e} and of the proof of $(iii)\Rightarrow (i)$ in
Theorem~\ref{thm:existence} which makes more precise the
equivalence of $(i)$ and $(ii)$ in Theorem~\ref{thm:existence}.

\begin{cor}
Let $n$ be a positive integer. 
A Seifert manifold whose base is not a sphere has a contact
structure $\xi$ with $t(\xi)=-n$ if and only if it has a
transverse contact structure $\xi'$ with $t(\xi')=-n$.
\end{cor}

In this section and the following one, we will use frequently the
idea of topological discretization which was first used in contact
geometry in \cite{Colin_chirurgie}. Recall that, for a surface
$F$, the pinched product $F\times_\partial[0,1]$ is obtained from
$F\times[0,1]$ by collapsing $\{x\}\times[0,1]$ for every
$x\in\partial F$. Two embedded surfaces in a 3--manifold are said
to be parallel if they bound a pinched product. Let $F$ be
embedded in $V$ and $\varphi$ be an isotopy of $V$ relative to
$\partial F$. Any time $t\in[0,1]$ has a neighborhood $J$ such
that $\varphi_{|J}$ moves all connected components of $F$ in
disjoint pinched products. Using this remark, one can show that
$\varphi$ is homotopic to a concatenation of isotopies
$\varphi^i$ which are relative to $\partial F$ and to all
connected components of the image of $F$ under the preceding
$\varphi^j$ but one denoted by $F_i$ and $\varphi^i_0(F_i)$ and
$\varphi^i_1(F_i)$ bound a pinched product.

If $V$ carries a contact structure $\xi$ then we can assume that
all the intermediate surfaces arising in the preceding process
are $\xi$--convex using the genericity of $\xi$--convex surfaces.

This discretization process will be called topological
discretization to avoid confusion with the more elaborate contact
discretization first used systematically in \cite{Honda_gluing} which
goes further by imposing---using
\cite{Giroux_transfos}[Lemma~15]---pinched products which are as
simple as possible from a contact point of view.

\begin{lem}
\label{lemma:calcul_e_disque}
  Let $V$ be a Seifert manifold with base an
orbifold disk and $W_1,\dots,W_r$ be fibered solid tori
such that $V\setminus(W_1\cup\dots\cup W_r)$ contains only
regular fibers.
 If $\xi$ is a contact structure which is tangent to the fibers
outside the $W_i$'s with twisting number $-n$ and $t(\xi_{|W_i})=-n$ for
every $i$ then $t(\xi)=-n$ (in particular $\xi$ is tight).
\end{lem}

\begin{proof}
Let $A_1,\dots,A_r$ be fibered annuli with boundary in $\partial V$ such 
that 
\[
V\setminus(A_1\cup\dots\cup A_r)=W_0'\cup\dots\cup W'_r
\]
where $W_i\subset W_i'$ for every $i\geq 1$.

\begin{figure}[ht]
 \labellist
 \small\hair 2pt
 \pinlabel $A_1$ at 25 44
 \pinlabel $A_2$ at 27 24
 \pinlabel $A_3$ at 48 28
 \pinlabel $N$ at 149 41
 \pinlabel $W$ at 149 21
 \endlabellist

  \begin{center}
    \includegraphics{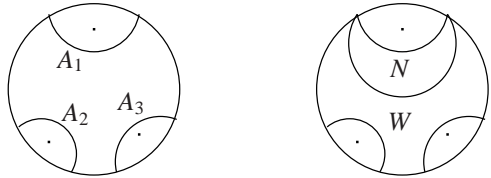}
  \end{center}
   \caption{Computing $t(\xi)$ over an orbifold disk}
  \label{fig:disk}
\end{figure}

Let $L_0$ be a fiber in $W_0'$, $L$ a vertical Legendrian curve and
$\varphi$ a $\partial$--isotopy such that $L=\varphi_1(L_0)$.

Using topological discretization, there exists a sequence of
annuli $A_i^j$, $0\leq j\leq K$ such that:
\begin{itemize}
  \item for every $i$, $A_i^0=A_i$, $A_i^K=\varphi_1(A_i)$;
  \item for every $j$ the $A_i^j$'s are disjoint $\xi$--convex
	annuli;
  \item for every $j$ there exists $i_j$ such that $A_i^{j+1}=A_i^{j}$ 
	  for all $i\neq i_j$ ;
  \item	$A_{i_j}^j$ and $A_{i_j}^{j+1}$ bound a pinched product.
\end{itemize}

At each step $j$ we equip $V$ with a Seifert fibration
$\partial$--isotopic to the original one such that all the $A_i^j$
are fibered.

We now prove by induction on $j$ that for every $j$ the
complement of the $A_i^j$ doesn't contain any vertical Legendrian
curve with $t>-n$ and each $A_i^j$ is divided by $2n$ curves 
going from one boundary component to the other one. This will prove 
the lemma since $L$ is in the complement of the $A_i^K$.

This statement is true at the initial step by assumption. Let's
assume that it holds at step $j$, set $A=A_{i_j}^j$,
$A'=A_{i_j}^{j+1}$ and let $N$ be the pinched product they bound. 
The annulus $A'$ is contained either in the
solid torus that $A$ splits from $V$ or in the solid torus
obtained from $V$ by removing all the solid tori split from $V$ 
by the $A_i^j$. Let's denote this solid torus containing $A'$ 
by $W$ and denote by $W'\subset W$ the complement of $N$ in $W$.

By the induction hypothesis, the dividing set of $A$ consists of
$2n$ curves going from one boundary component to the other
one. Since $A'$ is contained in $W$ we know by induction
hypothesis that it does not
contain any vertical Legendrian curve with twisting number
greater than $-n$ and since $\partial A'=\partial A$, 
Lemma~\ref{lemma:exist_ruling} guarantees that
the dividing set of $A'$ consists of $2n$ curves going from 
one boundary component of $A$ to the other one.

If $n>1$ then according to Proposition~\ref{prop:local_index}
the dividing set of $W'$ is isotopic to that of $W$ and the 
classification of tight contact structures on solid tori proves 
that there is a $\partial$--isotopy which brings $A$
to $A'$ through $\xi$--convex surfaces so nothing changes.

If $n=1$ then according to Proposition~\ref{prop:local_index} we only 
have to prove that $\xi$ restricted to all connected components of
the complement of the $A_i^{j+1}$'s is tight. All of these components 
but one are contained in solid tori which are known to be tight by 
the induction hypothesis.
The last component is in the union of a tight solid torus and a
tight toric annulus whose boundary components are divided by two
curves with homology classes $d=S+kF$ and $\dprime=S+k'F$
respectively so this component is tight according to 
Lemma~\ref{lemma:gluing_pour_calcul_e}. 
\end{proof}

\begin{proof}[Proof of Theorem~\ref{thm:e}]
Let $L$ be a vertical Legendrian curve in $V$. There exists a
finite cover of the base $B$ coming from a cover of the
underlying topological surface which induces a fibered covering
$\Vhat \to V$ such $L$ can be lifted to a Legendrian curve 
$\Lhat$ which is contained in a Seifert manifold with boundary $\Vhat'$ 
over an orbifold disk and $\Lhat$ is $\partial$--isotopic to
a regular fiber. Using Lemma~\ref{lemma:calcul_e_disque} we get
$t(L)=t(\Lhat)\leq -n$.
\end{proof}

We end this part with some remarks about the genus zero case. The
theorem above can be extended easily to spherical bases when
there is less than three exceptional fibers or when each $W_i$
has a unique tight contact structure, see 
Subsection~\ref{SS:tangentes}. Concerning small Seifert manifolds
with three exceptional fibers which are the subjects of many
recent papers, there are three cases: when $e_0(V)\geq 0$ then
according to Theorem~\ref{thm:existence} there exist no contact
structure with negative twisting number. When $e_0(V)\leq -3$
then condition (iii) of Theorem~\ref{thm:existence} can be met
only when $n=1$ and \cite{Wu_classif} (or
Theorem~\ref{thm:flexible}) combined with \cite{Wu_circles}
proves that the theorem is true in this case. The remaining cases
$e_0(V)\in\{-2,-1\}$ are much more difficult, very few results
are known especially when there are contact structures with
different negative twisting numbers.

As an example of what we can still get using our techniques,
let's consider the Brieskorn homology sphere $V=-\Sigma(2,3,6k-1)$
which has a Seifert structure with invariants
$(0,-2,(2,1),(3,2),(6k-1,5k-1))$ and has been studied for instance in
\cite{LM_Stein}. Our results prove that for every contact
structure $\xi$ on $V$, if $t(\xi)<0$ then 
\[t(\xi)\in\{-(6l+5),\; 0 \leq l \leq k-2\}.\]
There is exactly one isotopy class of
contact structures with $t=-(6(k-2)+5)$, its elements are universally
tight and isotopic to a tangent contact structure, see 
Subsection~\ref{SS:tangentes}. Theorem~\ref{thm:existence} gives at least one transverse
contact structure with $t\ge -5$ and $-5$ is the maximal possible
negative twisting number so there is at least one
isotopy class of universally tight contact structures with
$t=-5$. If the hypothesis $g>0$ can be removed from Theorems
\ref{thm:e} and \ref{thm:rigide} then the predicted number of
isotopy classes of tight contact structures with $t=-(6l+5)$ is
$k-l-1$ with two consisting of universally tight contact
structures when $t> -(6(k-2)+5)$. 
In any cases this is an upper bound on the number of
such isotopy classes.

Using Eliashberg--Gompf surgery \cite{Eliashberg_Stein,Gompf} 
and the slice Thurston--Bennequin inequality 
\cite{LM_Stein}[Corollary~4.2] it can be proved that the
predicted number is correct for $t=-5$ and that they are are all
Stein fillable and distinguished by their Stein fillings (I thank
Paolo Lisca for a very instructive conversation which led to this
result). A similar phenomenon occurs in a paper by J~Kim
\cite{Kim} where it is shown that the Seifert manifolds with
invariants $(0,-2,(2,1),(3,2),(6k-1,6k-3))$, $k\geq 2$ have 
$3k-5$ isotopy classes of tight contact structures. Using the
above arguments, we can show that all these contact structures 
have twisting number $-5$, as expected using the results of the
present paper. In this case the fact that $t=-5$ is the only
possible twisting number seems to be what makes the
classification feasible using such techniques.

It seems that there is no known counter-example to the following
statement: if a Stein fillable contact structure on a Seifert
manifold $V$ has twisting number $t_0<0$ then there is no contact
structure on $V$ with twisting number $t_0<t<0$. This can be
checked for instance on Lens spaces and on $T^3$ and we can note
that, according to P~Ghiggini in
\cite{Paolo_fillings}, ``Stein'' cannot be replaced by
``strongly'' in the above sentence.

\section{Classification}
\label{S:classification}

\subsection{Flexible case}

In this subsection we prove Theorem~\ref{thm:flexible} from the
introduction.

This result was announced by K~Honda for circle bundles ($r=0$) 
in \cite{Honda_II} 
but proved using the incorrect claim that these
contact structures can always be distinguished
by their homotopy classes as oriented plane fields.
Indeed, when $V$ has a contact structure which is tangent to the
fibers (i.e. when there exists $n$ such that $n e(V)=-\chi(B)$, see
\cite{Giroux_2001} or Theorem~\ref{thm:carac_tangentes}) there is
only one homotopy class of (oriented) plane fields transverse to
the fibers because the tangent contact structure is isotopic to
contact structures transverse to the fibers with either
orientation. If $-n<-1$ and $B$ has genus at least two then
$e<2g-2$. The mistake in \cite{Honda_II} arose from overlooking
the fact that fibers of $V$ have finite order in
$H_1(V,\Z)$ when $e\neq 0$.

\begin{lem}
Let $V$ be a Seifert manifold with $e_0(V)<2g-2$. Two
contact structures on $V$ with twisting number $-1$
always have normal forms with multi-index 
\[(2-2g-r,1,\dots,1)\] and the same $R$--class.
\end{lem}

\begin{proof}
We first prove that any such contact structure $\xi$ has a normal
form with multi-index $(2-2g-r,1,\dots,1)$. Let $\xi'$ be a
normal form of $\xi$ (obtained using
Proposition~\ref{prop:exist_normal_form}). According to
Proposition~\ref{prop:local_index}, the multi-index of $\xi'$
satisfies $x_i\leq 1$ for every $i\geq 1$. So, according to
Corollary~\ref{cor:exist_tore} and the flexibility lemma,  for
every $i\geq 1$, $W_i$ contains a torus $T_i$ parallel to the
boundary, ruled by vertical Legendrian curves with twisting
number $-1$ and divided by two curves whose homology class is
$S_i$. Consider a  vertical Legendrian curve $L_0$ outside the
$W_i$'s with $t(L_0)=-1$ and, for each $i$, a $\xi$--convex
annulus $A_i$ such that $\partial A_i$ is the union of $L_0$ and
vertical Legendrian curve on $T_i$. According to
Lemma~\ref{lemma:exist_ruling}, these annuli are
isotopic---relative to $L_0$ and the $T_i$'s---to annuli
intersecting only along $L_0$ and ruled by vertical Legendrian
curves with twisting number $-1$. One can then follow the proof
of the existence of normal forms to get a normal form with
$x_i=1$ for every $i\geq 1$. According to the Poincar\'e--Hopf
theorem, this normal form has $x_0=2-2g-r$.

The proof that if $\xi$ and $\xi'$ are in normal form with
this multi-index then $\xi$ has a normal form with 1--form
homotopic to that of $\xi'$ is exactly as in
\cite{Giroux_2001}[Lemma~3.8]
using the fact that $x_0<b$ so that $d_0=M_0+(x_0-b-1)F_0$ with
$x_0-b-1\leq -2$.
\end{proof}

Note that neither the above lemma nor Theorem~\ref{thm:flexible}
claims that every normal form with twisting number $-1$ has the
given multi-index---this would be false---but only that there is
a normal form with this multi-index in every isotopy class of
$\xi$ with $t(\xi)=-1$. This contrasts with Theorem~\ref{thm:rigide} where the multi-index is fixed by the Seifert
invariants and the twisting number.

\begin{proof}[Proof of Theorem~\ref{thm:flexible}]
  The preceding lemma proves that the number of contact
structures on $V$ with $t=-1$ is at most the number of contact
structures one can obtain by fixing a tangent contact structure
over $R$ with multi-index $(2-2g-r,1,\dots,1)$ and extending it
by a tight contact structure $\xi_i$ in each $W_i$. The proof
that this upper bound is the exact count follows a well-known
strategy---see eg \cite{Ghiggini_thesis,Wu_classif}---so we
only indicate the steps.

We use the Eliashberg--Gompf construction of Stein fillable
contact structures \cite{Eliashberg_Stein,Gompf} to
construct the right number of diffeomorphic Stein fillings of $V$
with different first Chern classes. The induced contact structures
are non-isotopic according to \cite{LM_Stein}[Corollary~4.2].
If $\xi$ is one of them then the Gompf diagram shows that 
$t(\xi)\geq -1$, and $t(\xi)<0$ according to the slice
Thurston--Bennequin inequality \cite{LM_Stein}[Theorem~3.4] so
$t(\xi)=-1$.

The count of transverse contact structures follows from this,
Proposition~\ref{prop:univ_tendue_positive} and 
Lemma~\ref{lemma:W_univ_tendue} using $d_0\cdot M_0>1$. 
Suppose that $\xi$ and $\xi'$ are in normal form with $t=-1$,
isotopic to transverse contact structures
and coincide on $R\times \S^1$ but are not isotopic.
The rotation of angle $\pi$ on $R\times \S^1$ pushes $\xi$ to
$-\xi'$ and normal forms in the $W_i$'s show that it 
can be extended to a diffeomorphism of $V$ still isotopic to the
identity and pushing $\xi$ to $-\xi'$.
\end{proof}

\subsection{Rigid case}

In this subsection we prove Theorem~\ref{thm:rigide} from the
introduction.

\begin{prop}
\label{prop:multi_index_can}
 Let $n$ be a positive integer and  $V$ be a Seifert manifold
whose base has genus $g$. 
If $e_0=2g-2$ or $n>1$ then every contact 
structure $\xi(\lambda,n,\xi_0,\dots,\xi_r)$
with $t(\xi_i)\leq -n$ for every $i$ has multi-index
$(nb,\lceil n\beta_1/\alpha_1\rceil,\dots,
\lceil n\beta_r/\alpha_r\rceil)$.
\end{prop}

\begin{proof}
If $n>1$ then the proposition follows directly from Proposition 
\ref{prop:local_index}. If $n=1$ then we need the additional
remark that the multi-index $(x_0,\dots,x_r)$ of $\xi$ 
satisfies $x_0-1<b$ and $x_i-1<\beta_i/\alpha_i$ if $i>0$ so 
$x_0\leq b$ and $x_i\leq 1$ if $i>0$ but also 
$\sum x_i=2-2g$ and $2-2g=b+r$ by hypothesis
so all the inequalities are equalities.
\end{proof}

The following proposition will be proved later in this
subsection. It is the only point where we need $g$ to be
positive. A proof in the $g=0$ case would lead to the extension
of Theorem~\ref{thm:rigide} to $g=0$.

\begin{prop}
\label{prop:unique_curve}
  Let  $n$ be a positive integer and $V$ a Seifert manifold 
whose base has genus $g>0$.
Assume that $n>1$ or $e_0(V)=2g-2$. Let $\xi$ be a contact
structure on $V$ with $t(\xi)=-n$. If $L_0$ and
$L_1$ are vertical Legendrian curves with $t(L_0)=t(L_1)=-n$
then every isotopy $L_t$ between $L_0$ and $L_1$ is 
homotopic to a Legendrian isotopy.
\end{prop}

We will also use the following 
flexibility lemma which is a special case of the general
flexibility lemma for families (\cite{Giroux_transfos}[Lemma~7]).

\begin{lem}
  \label{lemma:flexibility}
  Let $F$ be a closed surface in a contact 3--manifold
$(V,\xi)$. Suppose $\xi F$ is divided by a multi-curve $\Gamma$
and $\varphi$ is an isotopy such that $\varphi_t(F)$ is divided
by $\varphi_t(\Gamma)$ for every $t$. If 
$\xi(\varphi_1(F))=\varphi_1(\xi F)$ then $\varphi$ is homotopic to an
isotopy $\psi$ such that $\xi(\psi_t(F))=\psi_t(\xi F)$ for every $t$.
\end{lem}

\begin{proof}[Proof of Theorem~\ref{thm:rigide}]
If two contact structures $\xi$, $\xi'$ in normal forms with
$t=-n$ have the same $R$--class then $\xi'$ is isotopic through
contact structures in normal form to $\xi''$ with the same Euler
classes as $\xi'$ and $\lambda''=\lambda$. Because of
Proposition~\ref{prop:local_index}, each characteristic foliation
of a $\partial W_i$ is divided by two curves so $\xi''$ is
isotopic to $\xi$ according to
Theorem~\ref{thm:classif_tore_plein}~c).

Conversely, let $\xi$ and $\xi'$ be contact structures in normal
form on $V$ with twisting number $-n$ and suppose there is
an isotopy $\phi$ pushing $\xi'$ on $\xi$. We will simplify
$\phi$ in three steps. We make it relative to a fiber
$L_0\in R$ in step one, to a system of tori intersecting along
$L_0$ with a regular neighborhood isotopic to $R\times \S^1$ in
step two and to the $\partial W_i$'s in step three.

\step{1}

According to Proposition~\ref{prop:multi_index_can}, we can assume---up to an isotopy
of $\xi'$ among contact structures in normal form 
preserving its $R$--class and the Euler classes
of the $\xi'_i$'s---that $\xi'_{|\partial W_i}=\xi_{|\partial
W_i}$ for every $i$.

Let $L_0$ be a fiber in $R$. The fiber $L_0$ and its image $L_1$
by $\phi_1$ are $\xi$--Legendrian with maximal twisting number so
we can apply Proposition~\ref{prop:unique_curve} and the Legendrian isotopy
we get is induced by an isotopy $\psi$ such that $\psi_1=\phi_1$.

Since $\psi_t(L_0)$ is $\xi$--Legendrian for all $t$, there exists
an isotopy $\theta$ preserving $\xi$ which coincides with $\psi$
on $L_0$ for all time. Denote by $\phibar$ the isotopy given by
$\theta_t^{-1}\circ\psi_t$. This isotopy pushes $\xi'$ on $\xi$,
is relative to $L_0$, and we have
\[\xi\phibar_1(\partial W_i)=\xi\phi_1(\partial W_i)=\phi_1(\xi'\partial
W_i)=\phi_1(\xi\partial W_i)=\phibar_1(\xi\partial W_i).\]

\step{2}

Let $T_1,\dots,T_{2g+r}$ be fibered tori in $R\times \S^1$ 
intersecting along $L_0$ such that a regular neighborhood of
$T_1\cup\dots\cup T_{2g+r}$ is isotopic to $R\times \S^1$.

We now prove that $\phibar$ is homotopic to an isotopy
which moves the $T_i$'s through $\xi$--convex surfaces. The key is
Proposition~\ref{prop:multi_index_can} which gives
the multi-index independently of the normal form.

Using topological discretization, it is sufficient to prove that
if $\phibar$ is relative to all the $T_i$'s but one denoted by
$T$ and if $T':=\phibar(T)$ bounds a pinched product $N$ with $T$
then $\phibar$ is homotopic to an isotopy which, in addition,
moves $T$ through $\xi$--convex surfaces. Using the flexibility
lemma and Lemma~\ref{lemma:exist_ruling} we can assume that $T'$
is ruled by vertical Legendrian curves with twisting number $-n$.

The complement of $T_1\cup\dots\cup T_{2g+r}$ is the union of
$r+1$ solid tori. Let $W$ be the one that contains $T'$
and set $W'=W\setminus N$. The solid tori $W$ and $W'$ are
nested and they have isotopic dividing sets according to
Proposition~\ref{prop:multi_index_can}. The classification of
tight contact structures on solid tori then implies that
$\phibar$ is homotopic to an isotopy moving $T$ through convex
surfaces. 

This already proves that $\lambda$ is homotopic to $\lambda'$.
Moreover, because $\xi\phibar_1(T_i)=\phibar_1(\xi'T_i)$ the tori
$\phibar_1(T_i)$ are ruled by $\xi$--Legendrian curves so, 
using Lemma~\ref{lemma:flexibility}, we get an isotopy
$\varphi$ homotopic to $\phibar$ and such that 
$\xi\varphi_t(T_i)=\varphi_t(\xi T_i)$ for all $t$ and
every $i$. 

\step{3}

At all times, there is a regular neighborhood of 
$\bigcup\varphi_t(T_i)$ foliated by vertical Legendrian curves so
$\varphi$ is homotopic to 
$\widetilde\psi$ such that 
$\widetilde\psi_t(\partial W_i)$ is ruled by vertical curves for all $t$ and
every $i$. Using Lemma~\ref{lemma:flexibility}, we get an 
isotopy $\bar\psi$ homotopic to $\widetilde\psi$ such that 
\[\xi\bar\psi_t(\partial W_i)=\bar\psi_t(\xi\partial W_i)\]
for all $t$ and every $i$.

Let $\bar\theta$ be an isotopy preserving $\xi$ which coincides
with $\bar\psi$ on every $\partial W_i$, and denote by
$\widetilde\phi$ the isotopy given by
$\widetilde\phi_t=\bar\theta_t^{-1}\circ\bar\psi_t$.
This isotopy pushes $\xi'$ on $\xi$ and is relative to the
$\partial W_i$'s so that each $\xi_i$ is isotopic to the
corresponding $\xi'_i$ so they have the same relative Euler class.

The count of transverse contact structures goes as in the proof
of Theorem~\ref{thm:flexible} using
Proposition~\ref{prop:univ_tendue_positive},
Lemma~\ref{lemma:W_univ_tendue} and
Proposition~\ref{prop:local_index}. There is only one isotopy
class of transverse contact structure in each $R$--class if and
only if $d_i\cdot M_i=1$ for every $i$ and then this isotopy
class contains a tangent contact structure, see
Subsection~\ref{SS:tangentes}.
\end{proof}

\subsubsection*{Uniqueness of the maximally twisting vertical
curve}

In this paragraph we prove Proposition~\ref{prop:unique_curve}. 
As a first approach to this proposition, one can try to
consider a vertical $\xi$--convex torus containing $L_0$,
discretize the isotopy and prove that all bifurcations
can be eliminated. However this can not be straightforward since
one can always have bifurcations increasing the number of
dividing curves by folding inside an invariant neighborhood and
it is difficult to prove that a sequence of folding and
unfolding can be unraveled.
Here we use Ghiggini's trick introduced in
\cite{Ghiggini_thesis}[Proposition~5.4] and consider two tori
intersecting along $L_0$, discretize the motion of one of them
while constructing an isotopy of the second one through $\xi$--convex
surfaces. Here again, everything is based on the fact that the
involved Giroux cones contain only one integer point having the
right intersection with the fiber class.

\begin{lem}
  \label{lemma:anneau_traverse}
  Let $V$ be a Seifert manifold and $\xi$ be a contact structure
with negative twisting number $t(\xi)=-n$. Let $T$ and $F$ be
transverse vertical 
$\xi$--convex tori such that $T\cap F$ is a Legendrian
vertical curve which intersects the dividing set of $T$
efficiently.
If $n>1$ or $e_0(V)=2g-2$ then every isotopy relative to $T$
which sends $F$ to a $\xi$--convex torus $F'$ is homotopic 
to an isotopy relative to $T$ and moving $F$ through $\xi$--convex
surfaces.
\end{lem}

\begin{proof}
  Up to a change of Seifert structure by isotopy, we can assume
that $T$ is fibered. Let $-n'$ be the maximal twisting number
of Legendrian curves isotopic to the regular fibers relative to
$T$. Using topological discretization and the flexibility lemma,
we can assume that $F$ contains a vertical Legendrian curve $L$
with $t(L)=-n'$ and that $F$ and $F'$ bound a pinched product $N$
intersecting $T$ only along $L$. Let $\Sigma$ and $\Sigma'$ be
$\xi$--convex tori bounding regular neighborhoods of $T\cup F$ and
$T\cup N$ respectively, chosen so that they both contain a
vertical Legendrian curve $L_0$ with $t(L_0)=-n'$.

\vspace{3mm}

{\noindent\bf Claim:} the tori $\Sigma$ and $\Sigma'$ have isotopic
dividing sets with total homology class $2n'S-2F$.

\vspace{3mm}

Using this claim we can cut $V$ along $\Sigma'$, keep the regular
neighborhood of $T\cup N$ and fill it with a solid torus $W$ with
meridian class $S$ and a tight contact structure---there is no
choice here, up to isotopy---to get a contact manifold $V'$
diffeomorphic to $T^3$. Cutting $V'$ along $T$ we get a toric
annulus with a tight contact structure having the same dividing
set on both boundary components. Using
\cite{Giroux_2000}[Theorem~4.5] we see that $F$ and $F'$ are isotopic through
convex surfaces in $V'$ relative to $T\cup W$ so relative to the
boundary and the lemma is proved.

We now prove the claim. By maximality of $-n'$, we know that
$\Sigma$ and $\Sigma'$ have dividing sets with total homology
classes $2n'S+2kF$ and $2n'S+2k'F$ respectively. We will first
prove that $k=k'$ and then that $k=-1$. We consider two
cases:

If $n'>1$ then Corollary~\ref{cor:exist_tore} and Pick's formula
prove that $k\neq k'$ contradicts the
maximality of $-n'$ as in the proof of Proposition~\ref{prop:local_index}.

If $n'=1$ then $n=1$ and the proof is a variation of the proof of
Proposition~\ref{prop:multi_index_can}. We can construct a tori
system based at $L_0$---and not intersecting $\Sigma$ and
$\Sigma'$ anywhere else---which is ruled by vertical Legendrian
curves with twisting number $-1$. We equip $V$ with a Seifert
structure isotopic to the original one so that our tori system is
fibered. The tori $\Sigma$ and $\Sigma'$ intersect the
corresponding $W_0$ solid torus along annuli with common
boundaries. Using Proposition~\ref{prop:multi_index_can} we know
that the Giroux cone of $W_0$ contains only one integer point
with abscissa $1$ so that we can conclude using the
classification of tight contact structures on solid tori like in
the proof of Theorem~\ref{thm:rigide}.

So in both cases we proved that $k'=k$. It remains to prove that
$k=-1$. It is sufficient to construct a curve with homology class
$S$ which intersects only twice the dividing set of $\Sigma$
with the correct orientations. We know that the dividing set of
$T$ is made of curves traversing $T\setminus (T\cap F)$ and that 
$2n'$ dividing curves of $F\setminus (T\cap F)$ are
traversing. Our curve is constructed so that it intersects the dividing
curves only in the rounding regions. It starts just below (for the
fibers orientation) one of the traversing curves $C$ of $F$, follows
it then traverse $T$ without intersecting its dividing set then
traverse back $F$ below $C$ and traverse back $T$ before
closing up. There are two intersection points with the dividing
set of $\Sigma$ and the orientation is correct because all
traversing curves in the two copies of $F\setminus (T\cap F)$ are
oriented in the same way and because of the edge-rounding lemma
\cite{Honda_I}[Lemma~3.11].
\end{proof}

\begin{proof}[Proof of Proposition~\ref{prop:unique_curve}]
Let $\phi$ be an isotopy such that $L_t=\phi_t(L_0)$.
  Let $T$ and $F$ be $\xi$--convex vertical tori intersecting
transversely along 
$L_0$. Up to a modification of $\phi$ relative to $L_0$
we can assume that $T':=\phi_1(T)$ and $F':=\phi_1(F)$ are also
$\xi$--convex. By maximality of $t(L_0)$ (resp. $t(L_1)$) and the
flexibility lemma, $L_0$ (resp. $L_1$) intersects efficiently the 
dividing sets of $T$ and $F$ (resp. $T'$ and $F'$).

By topological discretization, $\phi$ is homotopic to a
concatenation of isotopies moving $T$ through
a sequence $T=T_0,T_1,\dots,T_N=T'$ of $\xi$--convex 
tori such that $T_i$ and $T_{i+1}$ bound a product.
Also, using the flexibility lemma at each step, we can assume that
every $T_i$ contains a vertical Legendrian curve $C_i$
intersecting efficiently the dividing set of $T_i$, $C_0=L_0$ and
$C_N=L_1$.

We now prove by induction that, for every $i$, there exists a torus
$F_i$ intersecting $T_i$ along $C_i$ which is isotopic to $F$
through $\xi$--convex surfaces.
For the initial step we take $F_0$ to be $F$. 

Once $F_i$ has been constructed we consider a torus $F_{i+1}$ which
is isotopic to $F_i$ relative to $T_i$ and intersects $T_{i+1}$
along $C_{i+1}$. According to Lemma~\ref{lemma:anneau_traverse},  
$F_{i+1}$ is isotopic to $F_i$ through $\xi$--convex surfaces.

We apply Lemma~\ref{lemma:anneau_traverse} one last time to get
an isotopy between $F_N$ and $F'$ and the proposition is proved.
\end{proof}

\begin{rem}
About Legendrian knots:

 \begin{itemize}
  \item Lemma~\ref{lemma:anneau_traverse} can also be used 
	to prove that any vertical Legendrian curve is a
	stabilization of one with maximal twisting number so we
	have a complete classification of vertical Legendrian knots
	under the hypotheses of Theorem~\ref{thm:rigide}. 
  \item When $t(\xi)=-1$ and $e_0=2g-2$, Theorem~\ref{thm:flexible} 
	can be used to prove that if $t(L_0)=t(L_1)=-1$	
	then there exists a contactomorphism which sends $L_0$ to
	$L_1$ and is isotopic to the identity---although possibly not
	through contactomorphisms. If one could prove the existence
	of a Legendrian isotopy bringing $L_0$ to $L_1$ then it would
	be possible to prove Theorem~\ref{thm:flexible} in the spirit
	of the proof of Theorem~\ref{thm:e} (without using
	Seiberg--Witten theory).
 \end{itemize}
\end{rem}

\subsection{Tangent contact structures}
\label{SS:tangentes}

Among transverse contact structures on a Seifert manifold are
those which are isotopic to tangent contact structures 
such as the standard contact structure on $\S^3$.
Conversely, any tangent contact structure can be perturbed by a
$C^\infty$--small isotopy to be positively or negatively
transverse. The following theorem explain their special role
among contact structures with negative twisting number.
It is a direct consequence of the previous results.

\begin{thm}
  \label{thm:carac_tangentes}
  Let $V$ be a Seifert manifold with invariants
$$(g,b,(\alpha_1,\beta_1),\dots,(\alpha_r,\beta_r))$$ 
and $n$ be a positive integer. The following conditions are
equivalent:
\begin{enumerate}
  \item $ne(V)=-\chi(B)$ and $n\beta_i\equiv 1 \mod \alpha_i$ for
	all $i$;
  \item $V$ carries a tangent contact structure
	with twisting number $-n$.
\end{enumerate}

In addition, if 1 (or 2) is satisfied then every
contact structure on $V$ with twisting number
$-n$ is isotopic to a tangent one.

Also, if $V$ carries a contact structure with twisting
number $-n$ and $n>\max(\alpha_i)$ then 1 and 2 are satisfied and
this contact structure is isotopic to a tangent one.
\end{thm}

\begin{proof}
 Suppose that $V$ has a tangent contact structure $\xi$  
and $t(\xi)=-n$. Up to isotopy among tangent contact structures,
this $\xi$ is in normal form and, according to
Lemma~\ref{lemma:W_univ_tendue}, 
the multi-index $(x_0,\dots,x_r)$ of $\xi$ satisfies
$1=d_i\cdot M_i=n\beta_i-\alpha_i(x_i-1)$ so that
$n\beta_i \equiv 1 \mod \alpha_i$ for every $i$. Also
\begin{eqnarray*}
  ne(V) &=& - \sum_{i=0}^r\frac{n\beta_i}{\alpha_i} 
  = -\sum_{i=0}^r \left(\frac{1}{\alpha_i}+(x_i-1)\right) \\
  &=& -\left(\sum_{i=0}^r x_i - (r+1) + \sum_{i=0}^r
  \frac{1}{\alpha_i} \right) 
  = -\chi(B).
\end{eqnarray*}

Conversely if these conditions are met and $\xi$ is in normal
form we consider two cases.
If $n>1$ then Proposition~\ref{prop:local_index} and 
Lemma~\ref{lemma:W_univ_tendue} show that each $\xi_i$ is isotopic to a
tangent contact structure.
If $n=1$ then $\beta_i=1$ for every $i\geq 1$ because
$\beta_i\equiv 1 \mod \alpha_i$ and this combines with
$e(V)=-\chi(B)$ to give $-b= 2g-2+r$ so 
$e_0(V)=2g-2$ and according to Proposition~\ref{prop:multi_index_can} 
we can apply again Lemma~\ref{lemma:W_univ_tendue}.

Note that this proof of  $1)\Rightarrow 2)$ gives also the second assertion.

The last assertion of the theorem follows from Proposition 
\ref{prop:local_index} and Lemma~\ref{lemma:W_univ_tendue}.
\end{proof}

These contact structures are classified up to isotopy by Theorem 
\ref{thm:rigide} since, according to the preceding theorem, if
$\xi$ is of tangent type and $t(\xi)=-1$ then $e(V)=-\chi(B)$
and $\beta_i=1$ for every $i\geq 1$ so $e_0(V)=2g-2$.

It is plausible that the isotopy classes containing tangent contact 
structures are exactly the ones containing both positively and
negatively transverse contact structures. This would follow from the fact
that a positively transverse contact structure has a normal form
with positively transverse $\xi_i$'s but this is not what
Propositions \ref{prop:exist_normal_form} and
\ref{prop:univ_tendue_positive} give us.

Tangent contact structures also have the virtue that their
twisting number is easy to compute: it is given by the twisting
number of regular fibers, see \cite{Giroux_2001}[Lemma~3.6] 
which---using a cover by a circle bundle---only leaves the case of
certain Lens spaces which can be dealt with in the spirit of the
proof of Proposition~\ref{prop:univ_tendue_positive}. This remark
and the fact expressed in the above theorem that all contact
structures with sufficiently low twisting number $-n$ on a given
Seifert manifold satisfies $ne(V)=-\chi(B)$ can be used to get
a list of eight Seifert manifolds---with $\chi(B)=0$---which are
exactly the Seifert manifolds having an infinite family of
contact structures with distinct negative twisting numbers. They
are the cooriented contact elements bundles of the parabolic
orbifolds with either orientation (for $T^2$ and the pillowcase
the two orientations give the same Seifert manifold).

A better understanding of tangent contact structures and an
alternative proof of the equivalence of 1) and 2) in the previous
theorem come from the two following statements which can be proved 
exactly as in the circle bundle case using equivariance of all the
constructions. The second one is a slight generalization of 
\cite{Giroux_2001}[Proposition~3.3] which actually was first
discussed by R~Lutz in \cite{Lutz_pivot}.

\begin{lem}
\label{lemma:xi_B}
  Let $B$ be a 2--dimensional orbifold. The cooriented contact
elements bundle $\STB$ has a canonical contact structure $\xi_B$
which is tangent to the fibers with twisting number $-1$ and is
invariant under any fibered diffeomorphism lifted from $B$.
\end{lem}

\begin{prop}
\label{prop:rev_giroux_general}
 Let $V$ be a Seifert manifold with base $B$. The map which
associates to each covering map $\rho\co V \to \STB$ fibered over
the identity the contact structure $\rho^*\xi_B$ is a bijection
onto the space of tangent contact structures.
\end{prop}

%
%
%
\bibliographystyle{gtart}

\bibliography{gcs}
\end{document}